\documentclass{article}
\usepackage{amsmath}
\usepackage{arxiv}

\usepackage[utf8]{inputenc} 
\usepackage[T1]{fontenc}    
\usepackage{hyperref}       
\usepackage{url}            
\usepackage{booktabs}       
\usepackage{amsfonts}       
\usepackage{nicefrac}       
\usepackage{microtype}      
\usepackage{cleveref}       
\usepackage{lipsum}         
\usepackage{graphicx}
\usepackage{doi}

\usepackage{tikz}
\usepackage{subcaption}
\usepackage{graphicx}
\usepackage{epstopdf}
\usepackage{algorithmic}

\title{Well-balanced Central Scheme for the System of MHD Equations with Gravitational source term}

\author{ Farah Kanbar \\
University of Wuerzburg\\
Wuerzburg, Germany\\
 \AND
 Rony Touma \\
 Lebanese American University \\
 Beirut, Lebanon \\
  \AND
 Christian Klingenberg \\
 University of Wuerzburg \\
 Wuerzburg, Germany \\
}

\renewcommand{\shorttitle}{Well-balanced scheme for MHD}

\hypersetup{
pdftitle={Well-balanced Central Scheme for the System of MHD Equations with Gravitational source term},
pdfsubject={q-bio.NC, q-bio.QM},
pdfauthor={Farah Kanbar, Rony Touma, Christian Klingenberg},
pdfkeywords={First keyword, Second keyword, More},
}

\begin{document}
\maketitle

\begin{abstract}
  A well-balanced second order finite volume central scheme for the magnetohydrodynamic (MHD) equations with gravitational source term is developed in this paper. The scheme is an unstaggered central scheme that evolves the numerical solution on a single grid and avoids solving Riemann problems at the cell interfaces using ghost staggered cells. A subtraction technique is used on the conservative variables with the support of a known steady state in order to manifest the well-balanced property of the scheme. The divergence-free constraint of the magnetic field is satisfied after applying the constrained transport method (CTM) for unstaggered central schemes at the end of each time-step by correcting the components of the magnetic field. The robustness of the proposed scheme is verified on a list of numerical test cases from the literature.
\end{abstract}


\section{Introduction}
Ideal Magnetohydrodynamics (MHD) equations model problems in physics and astrophysics. The MHD system is a combination of the Navier-Stokes equations of fluid dynamics and the Maxwell equations of electromagnetism. A gravitational source term is added to the ideal MHD equations in two space dimensions in order to model more complicated problems arising in astrophysics and solar physics such as modeling wave propagation in idealized stellar atmospheres \cite{Retal,Tetal}.  
From electromagnetic theory, the magnetic field $\textbf{B}$ must be solenoidal i.e. $\nabla\cdot \textbf{B}=0$ at all times. The divergence-free constraint on the magnetic field reflects the fact that magnetic mono-poles have not been observed in nature. The induction equation for updating the magnetic field imposes the divergence on the magnetic field. Hence, a numerical scheme for the MHD equations should maintain the divergence-free property of the discrete magnetic field at each time-step. Numerical schemes usually fail to satisfy the divergence-free constraint and numerical instabilities and unphysical oscillations may be observed \cite{Toth}. Several methods were developed to overcome this issue. The projection method, in which the magnetic field is projected into a zero divergence field by solving an elliptic equation at each time step \cite{Brackbill}.\\
Another procedure is the Godunov-Powell procedure \cite{Powell2, Powell3, Fuchs2010}, where the Godunov-Powell form of the system of the MHD equations is discretized instead of the original system. The Godunov-Powell system has the divergence of the magnetic field as a part of the source term. Hence, divergence errors are transported out of the domain with the flow.\\
A third approach is the constrained transport method (CTM) \cite{CTM1,CTM2,CTM3}. The CTM was modified from its original form to the case of staggered central schemes \cite{ArminjonTouma2005}. It was later extended to the case of unstaggered central schemes \cite{Touma2010}. Hence, a numerical scheme for the MHD equations should maintain the divergence-free property of the discrete magnetic field at each time-step. 
A finite volume second-order accurate unstaggered central scheme is used to model the MHD equations with a gravitational source term. Finite volume central schemes were first introduced in 1990 by Nessyahu and Tadmor (NT) \cite{ARTICLE3_NT1}. The NT scheme is based on evolving  piecewise linear numerical solution on two staggered grids. The most significant property of central schemes is that they avoid solving Riemann problems arising at the cell interfaces. Our scheme is unstaggered central (UC) type scheme that was first developed in \cite{jian_al,amc09}. These schemes allow the evolution of the numerical solution on a single grid instead of using two different grids. UC schemes were first developed for hyperbolic systems of conservation laws and then extended to hyperbolic systems of balance laws\cite{ToumaKoleyKlingenberg,amc2012,touma2015well,APNUM2010}. The UC schemes introduced the possibility of avoiding solving Riemann problems and switching between two grids. The approach is achieved by the help of ghost staggered cells used implicitly to avoid Riemann problems at the cell interfaces.\\
In the presence of a gravitational source term on the right hand side of the MHD system, one has to consider a well-balanced technique that provides the numerical scheme with the ability to preserve hydrostatic equilibrium. In this paper we extend the reconstruction technique on the conservative variables, previously developed in  \cite{jonas,APNUM2020} for the system of Euler equations, for the system of MHD equations. The idea is to  evolve the error function between the vector of conserved variables and a given steady state, instead of evolving the vector of conserved variables. This error function is defined as $\Delta \textbf{U}=\textbf{U}-\tilde{\textbf{U}}$ , where $\tilde{\textbf{U}}$ is a given steady state. Knowing the steady state (analytically or numerically) is a key ingredient for the implementation of the proposed scheme.\\
The paper is divided into the following sections. The MHD model is presented in section \ref{sec::the model} and the finite volume scheme is described in section \ref{sec::the scheme} followed by the CTM in section \ref{sec::CTM}. Numerical experiments are illustrated in section \ref{sec::numerical results} and finally some concluding remarks and future work are given in section \ref{sec::conclusion}.

\section{The model}
\label{sec::the model}
The system of MHD  equations with gravitational source term in two space dimensions is given by:
\begin{equation}\label{MHDBL}
	\begin{cases}
		\textbf{U}_t+F(\textbf{U})_x+G(\textbf{U})_y=S(\textbf{U}), \hspace{0.5 cm}(x,y)\in\Omega\subset\mathbb{R}^2,\hspace{0.1 cm}t>0.\\
		\textbf{U}(x,y,0)=\textbf{U}_0(x,y),
	\end{cases}
\end{equation}
where
$$\textbf{U}= \left(  
\begin{array}{ c c }
	\rho \\
	\rho u_1 \\
	\rho u_2 \\
	\rho u_3 \\
	E \\
	B_1 \\
	B_2 \\
	B_3\\
	
\end{array} \right),~
F(\textbf{U})=\left(  
\begin{array}{ c c }
	\rho u_1 \\
	\rho u_1^2+\Pi_{11}\\
	\rho u_1 u_2+\Pi_{12}\\
	\rho u_1 u_3+\Pi_{13}\\
	Eu_1+u_1\Pi_{11}+u_2\Pi_{12}+u_3\Pi_{13} \\
	0 \\
	\Lambda_2\\
	-\Lambda_3\\
\end{array} \right),$$
$$
G(\textbf{U})=\left(  
\begin{array}{ c c }
	\rho u_2\\
	\rho u_2u_1+\Pi_{21}\\
	\rho u_2^2+\Pi_{22}\\
	\rho u_2 u_3+\Pi_{23}\\
	Eu_2+u_1\Pi_{21}+u_2\Pi_{22}+u_3\Pi_{23}  \\
	-\Lambda_3\\
	0 \\
	\Lambda_1\\
\end{array} \right), 
S(\textbf{U})=\left(  
\begin{array}{ c c }
	0 \\
	0 \\
	-\rho \phi_{y} \\
	0\\
	-\rho u_2 \phi_{y}  \\
	0 \\
	0\\
	0\\
\end{array}\right). $$
Here $\rho$ is the fluid density, $\rho \textbf{u}$ is the momentum with $\textbf{u}=(u_1,u_2,u_3)$, $p$ is the pressure, $\textbf{B}=(B_1,B_2,B_3)$ is the magnetic field, and $E$ is the kinetic and internal energy of the fluid given by the following equation $E=\frac{p}{\gamma-1}+\frac{1}{2}\rho |\textbf{u}|^2+\frac{1}{2}|\textbf{B}|^2$ with $\gamma$ the ratio of specific heats. $\phi=\phi(x,y)$, with $\phi_{x}=0$ and $\phi_{y}=g$, is the gravitational potential and it  is a given function. 
The conservation of the total energy (internal, kinetic and magnetic) has the gravitational potential energy as a source term.
$\Lambda=\textbf{u} \times \textbf{B}$, $\Pi_{11}, \Pi_{22}$ and $\Pi_{33}$ are the diagonal elements of the total pressure tensor and $\Pi_{12}, \Pi_{13}$ and $\Pi_{23}$ are the off-diagonal tensor are given by the following formulas:\\
$\Pi_{ii}=p+\frac{1}{2}(B_j^2+B_k^2-B_i^2)$ and $\Pi_{ij}=-\frac{1}{2}B_{i}B_{j}$, for $i,j,k=1, 2, 3.$\\
To determine the time-step using the CFL condition, we present the eigenvalues of the flux jacobian in the $x$-direction,\\
$\lambda_1=u_1-c_f$, $\lambda_2=u_1-b_1$, $\lambda_3=u_1-c_s$, $\lambda_4=u_1$, $\lambda_5=u_1$, $\lambda_6=u_1+c_s$, $\lambda_7=u_1+b_1$, $\lambda_8=u_1+c_f$. The eigenvalues of the flux jacobian in the $y$-direction are analogously defined.\\
Here, 
\begin{equation}
	c_f=\sqrt{\frac{1}{2}\left(a^2+b^2+\sqrt{\left(a^2+b^2\right)^2-4a^2b_1^2}\right)},
\end{equation}
and 
\begin{equation}
	c_s=\sqrt{\frac{1}{2}\left(a^2+b^2-\sqrt{\left(a^2+b^2\right)^2-4a^2b_1^2}\right)},
\end{equation}
are respectively the fast and slow wave speeds with $a=\sqrt{\frac{\gamma p}{\rho}}$ is the sound speed and $b=\sqrt{b_1^2+b_2^2+b_3^2}$ with $b_i=\frac{B_i}{\sqrt{\rho}}, i \in \{  1,2,3 \}.$ For additional reading on the hyperbolic analysis of the system, readers are refered to \cite{Godunov,Powell1}.
\section{The unstaggered two-dimensional finite volume central scheme}
\label{sec::the scheme}
\begin{figure}
	\centering
	\begin{tikzpicture}
		\fill [green!20] (0,0) rectangle (2,2);
		\fill [blue!20] (1,1) rectangle (3,3);
		\draw[step=1cm,black,thin,dashed] (0,0) grid (4,4);
		
		\fill (2,2)[blue] circle [radius=2pt];
		\fill (1,1)[green] circle [radius=2pt];
		
		\node at (1,0.6) {$(x_{i-\frac{1}{2}},y_{j-\frac{1}{2}})$};
		\node at (2,1.5)  {$(x_{i},y_{j})$};
		\node at (-1,4) {$(x_{i-1},y_{j+1})$};
		\node at (5,4) {$(x_{i+1},y_{j+1})$};
		\node at (5,0) {$(x_{i+1},y_{j-1})$};
		\node at (-1,0) {$(x_{i-1},y_{j-1})$};
	\end{tikzpicture}
	\caption{The cells of the main grid $C_{i,j}$ (blue cell) and of the staggered grid $D_{i-\frac{1}{2},j-\frac{1}{2}}$(green cell).}
	\label{2Dgrid}
\end{figure}
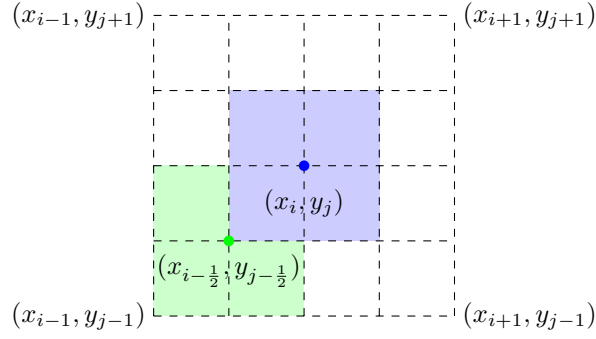

\begin{figure}[ht]
	\centering
	\begin{tikzpicture}
		\draw (0,0) -- (3,0) -- (3,3) -- (0,3) -- (0,0);
		\node at (0,-0.5) {$(x_{i},y_{j})$};
		\fill (0,0)[black] circle [radius=2pt];
		\node at (3,-0.5) {$(x_{i+1},y_{j})$};
		\fill (3,0)[black] circle [radius=2pt];
		\node at (3.5,3.5) {$(x_{i+1},y_{j+1})$};
		\fill (3,3)[black] circle [radius=2pt];
		\node at (-0.5,3.5) {$(x_{i},y_{j+1})$};
		\fill (0,3)[black] circle [radius=2pt];
		\draw [blue, -> ] (1.5,3) -- (1.5,4) node [right] {$\textbf{n}=(0,1)$};
		\draw [blue, -> ] (1.5,0) -- (1.5,-1)node [right] {$\textbf{n}=(0,-1)$};
		\draw [blue, -> ](3,1.5) -- (4,1.5)node [right]{$\textbf{n}=(1,0)$};
		\draw [blue, -> ] (0,1.5) -- (-1,1.5)node [left]{$\textbf{n}=(-1,0)$};
	\end{tikzpicture}
	\caption{The boundary $\partial R_{xy}$ and the outward pointing unit normal vector $\textbf{n}=(n_x,n_y)$ on each side of the boundary.}
	\label{BD2D}
\end{figure}
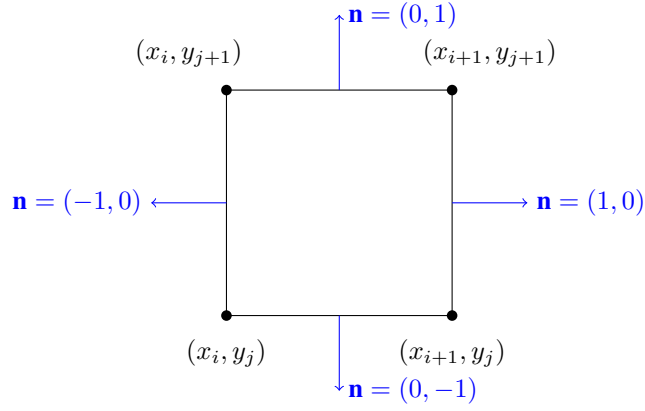

We consider a Cartesian decomposition of the computational domain $\Omega$ where the control cells are the rectangles 
$C_{i,j}=\left[x_{i-\frac{1}{2}},x_{i+\frac{1}{2}}\right]\times\left[y_{j-\frac{1}{2}},y_{j+\frac{1}{2}}\right]$ centered at the nodes $(x_i,y_j)$. We define the dual staggered cells $D_{i+\frac{1}{2},j+\frac{1}{2}}=\left[x_i,x_{i+1}\right]\times\left[y_j,y_{j+1}\right]$ centered at $(x _{i+\frac{1}{2}},y _{j+\frac{1}{2}})$.
Here, $x _{i+\frac{1}{2}}=x_i+\frac{\Delta x}{2} $ and $y_{j+\frac{1}{2}}=y_j+\frac{\Delta y}{2}$, where $\Delta x=x_{i+\frac{1}{2}}-x_{i-\frac{1}{2}}$ and $\Delta y=y_{j+\frac{1}{2}}-y_{j-\frac{1}{2}}$. The visualization of the 2D grids is given in figure \ref{2Dgrid}.
Before proceeding with the derivation of the 2D numerical method, and for convenience, we introduce the average value notations:
$$\overline{\rho}_{i,j+\frac{1}{2}}=\frac{\rho_{i,j}+\rho_{i,j+1}}{2},
\overline{\rho}_{i+\frac{1}{2},j}=\frac{\rho_{i,j}+\rho_{i+1,j}}{2}, \overline{\rho}_{i,(j)}=\frac{\rho_{i,j+\frac{1}{2}}+\rho_{i,j-\frac{1}{2}}}{2}$$
$$\overline{\rho}_{(i),j}=\frac{\rho_{i+\frac{1}{2},j}+\rho_{i-\frac{1}{2},j}}{2},\quad
[[\rho]]_{i,j+\frac{1}{2}}=\rho_{i,j+1}-\rho_{i,j}$$
$$[[\rho]]_{i+\frac{1}{2},j}=\rho_{i+1,j}-\rho_{i,j},~
[[\rho]]_{i,(j)}=\rho_{i,j+\frac{1}{2}}-\rho_{i,j-\frac{1}{2}}, ~
[[\rho]]_{(i),j}=\rho_{i+\frac{1}{2},j}-\rho_{i-\frac{1}{2},j}.$$
We assume that $\tilde{\textbf{U}}$ is a given stationary solution of system \eqref{MHDBL} and we define $\Delta\textbf{U}=\textbf{U}-\tilde{\textbf{U}}$.
We substitute $\textbf{U}=\Delta\textbf{U}+\tilde{\textbf{U}}$ in the balance law \eqref{MHDBL}, we obtain:
\begin{equation}
	(\Delta\textbf{U})_t+F(\Delta\textbf{U}+\tilde{\textbf{U}})_x+G(\Delta\textbf{U}+\tilde{\textbf{U}})_y=S(\Delta\textbf{U}+\tilde{\textbf{U}},x,y).
	\label{BLdelta2D}
\end{equation}
On the other hand, since $\tilde{\textbf{U}}$ is a stationary solution, then balance law in \eqref{MHDBL} reduces to
\begin{equation}
	F(\tilde{\textbf{U}})_x+G(\tilde{\textbf{U}})_y=S(\tilde{\textbf{U}},x,y).
	\label{BLtild2D}
\end{equation}
Subtracting equation \eqref{BLtild2D} from equation \eqref{BLdelta2D}, we obtain
\begin{multline}\label{BLdelta2S2D}
	(\Delta\textbf{U})_t+[F(\Delta\textbf{U}+\tilde{\textbf{U}})-F(\tilde{\textbf{U}})]_x+[G(\Delta\textbf{U}+\tilde{\textbf{U}})-G(\tilde{\textbf{U}})]_y\\
	=S(\Delta\textbf{U}+\tilde{\textbf{U}},x,y)-S(\tilde{\textbf{U}},x,y).
\end{multline}
Using the fact that the source term $S(\textbf{U},x,y)$ in \eqref{MHDBL} is linear in terms of the conserved variables, then equation \eqref{BLdelta2S2D} reduces to
\begin{align}
	(\Delta\textbf{U})_t+[F(\Delta\textbf{U}+\tilde{\textbf{U}})-F(\tilde{\textbf{U}})]_x+[G(\Delta\textbf{U}+\tilde{\textbf{U}})-G(\tilde{\textbf{U}})]_y=S(\Delta\textbf{U},x,y).
	\label{BLS2D}
\end{align}
The proposed numerical scheme consists of evolving the balance law \eqref{BLS2D} instead of evolving the balance law in system \eqref{MHDBL}. \\
The numerical solution $\textbf{U}$ will be then obtained using the formula $\textbf{U}=\Delta\textbf{U}+\tilde{\textbf{U}}$. The numerical scheme that we shall use to evolve $\Delta\textbf{U}(x,y,t)$ follows a classical finite volume approach; it evolves a piecewise linear function $\mathcal{L}_{i,j}(x,y,t)$ defined on the control cells $C_{i,j}$ and used to approximate the analytic solution $\Delta\textbf{U}(x,y,t)$ of system \eqref{MHDBL}. Without any loss of generality we can assume that $\Delta\textbf{U}_{i,j}^{n}$ is known at time $t^n$ and we define $\mathcal{L}_{i,j}(x,y,t^{n})$ on the cells $ C_{i,j}$ as follows.
\begin{align*}
	\mathcal{L}_{i,j}(x,y,t^{n})=\Delta\textbf{U}_{i,j}^{n}+(x-x_i)\frac{(\Delta\textbf{U}_{i,j}^{n,x})'}{\Delta x}+(y-y_j)\frac{(\Delta\textbf{U}_{i,j}^{n,y})'}{\Delta y}, \hspace{0.5 cm} \forall (x,y) \in C_{i,j},
\end{align*}
where $\frac{(\Delta\textbf{U}_{i,j}^{n,x})'}{\Delta x}$ and $\frac{(\Delta\textbf{U}_{i,j}^{n,y})'}{\Delta y}$ are limited numerical gradients approximating\\ $\frac{\partial\Delta\textbf{U}}{\partial x}(x,y_j,t^n)|_{x=x_i}$  and $\frac{\partial\Delta\textbf{U}}{\partial y}(x_i,y,t^n)|_{y=y_j}$, respectively, at the point $ (x_i, y_{j}, t^{n})$. 
In order to approximate the spatial numerical derivatives, the (MC-$\theta$) limiter is considered which is defined as 
\begin{align}\label{MC}
	(\Delta\textbf{u}_{i}^{n})'=\text{minmod}\left[\theta\left(\Delta\textbf{u}_{i}^{n}-\Delta\textbf{u}_{i-1}^{n}\right),
	\frac{\Delta\textbf{u}_{i+1}^{n}-\Delta\textbf{u}_{i-1}^{n}}{2},\theta\left( \Delta\textbf{u}_{i+1}^{n}-\Delta\textbf{u}_{i}^{n}\right)\right]
\end{align}
where $\theta$ is a parameter such that $1<\theta <2$, while the minmod function is defined as:
\[ \text{minmod}(a,b,c)=\begin{cases}
	\text{sign}(a)\text{min}\{|a|,|b|,|c|\}, \hspace{0.3 cm} \text{if} \hspace{0.1 cm} \text{sign}(a)=\text{sign}(b)=\text{sign}(c)\\
	0, \hspace{0.5 cm } \text{Otherwise.}
\end{cases}
\]
The (MC-$\theta$) limiter \eqref{MC} is used to compute the quantities $(\Delta\textbf{U}_{i,j}^{n,x})'$ and $(\Delta\textbf{U}_{i,j}^{n,y})'$  in order to avoid spurious oscillations. Next, we integrate the balance law \eqref{BLS2D} over the rectangular box $R_{i+\frac{1}{2},j+\frac{1}{2}}^n=D_{i+\frac{1}{2},j+\frac{1}{2}}\times[t^n,t^{n+1}]$,
\begin{multline}\label{IntBL2Da}
	\iiint_{R_{i+\frac{1}{2},j+\frac{1}{2}}}(\Delta\textbf{U})_t+[F(\Delta\textbf{U}+\tilde{\textbf{U}})-F(\tilde{\textbf{U}})]_x+[G(\Delta\textbf{U}+\tilde{\textbf{U}})-G(\tilde{\textbf{U}})]_y dR\\
	=\iiint_{R_{i+\frac{1}{2},j+\frac{1}{2}}}S(\Delta\textbf{U},x,y)dR.
\end{multline}
We use the fact that $\Delta\textbf{U}$ is approximated using piecewise linear interpolants similar to $\mathcal{L}_{i,j}$ on the cells $C_{i,j}$; following the derivation of the unstaggered central schemes in \cite{amc09}, equation \eqref{IntBL2Da} is rewritten as:
\begin{multline}\label{IntBL2Db}
	\Delta\textbf{U}_{i+\frac{1}{2},j+\frac{1}{2}}^{n+1}=\Delta\textbf{U}_{i+\frac{1}{2},j+\frac{1}{2}}^{n}-\frac{1}{\Delta x \Delta y}\iiint_{R_{i+\frac{1}{2},j+\frac{1}{2}}}[F(\Delta\textbf{U}+\tilde{\textbf{U}})-F(\tilde{\textbf{U}})]_x\\
	+[G(\Delta\textbf{U}+\tilde{\textbf{U}})-G(\tilde{\textbf{U}})]_y dR+\frac{1}{\Delta x \Delta y} \iiint_{R_{i+\frac{1}{2},j+\frac{1}{2}}}S(\Delta\textbf{U},x,y)dR.
\end{multline}
For the flux integrals, we apply the divergence theorem that converts the volume integral into a surface integral. Equation \eqref{IntBL2Db} becomes then:
\begin{multline}\label{IntBL2Dc}	\Delta\textbf{U}_{i+\frac{1}{2},j+\frac{1}{2}}^{n+1}=\Delta\textbf{U}_{i+\frac{1}{2},j+\frac{1}{2}}^{n}-\frac{1}{\Delta x \Delta y}\int_{t^n}^{t^{n+1}} \int_{\partial R_{xy} }[F(\Delta\textbf{U}+\tilde{\textbf{U}})-F(\tilde{\textbf{U}})]\cdot n_x  dAdt\\
	-\frac{1}{\Delta x \Delta y}\int_{t^n}^{t^{n+1}} \int_{\partial R_{xy}} [G(\Delta\textbf{U}+\tilde{\textbf{U}})-G(\tilde{\textbf{U}})]\cdot n_y  dAdt\\
	+\frac{1}{\Delta x \Delta y}\iiint_{R_{i+\frac{1}{2},j+\frac{1}{2}}}S(\Delta\textbf{U},x,y)dR,
\end{multline}
where $R_{xy}=\left[x_i,x_{i+1}\right]\times\left[y_i,y_{i+1}\right]$, and $\textbf{n}=(n_x,n_y)$ is the outward pointing unit normal at each point on the boundary $\partial R_{xy}$(the boundary of $R_{xy}$), see figure \ref{BD2D}.
The integral of the source term is being approximated using the midpoint quadrature rule both in time and space:
\begin{align}
	\iiint_{R_{i+\frac{1}{2},j+\frac{1}{2}}}S(\Delta\textbf{U})dR=\Delta x \Delta y \Delta tS(\Delta\textbf{U}_{i,j}^{n+\frac{1}{2}},\Delta\textbf{U}_{i+1,j}^{n+\frac{1}{2}},\Delta\textbf{U}_{i,j+1}^{n+\frac{1}{2}},\Delta\textbf{U}_{i+1,j+1}^{n+\frac{1}{2}}),
\end{align}
with
\begin{multline*}
	S(\Delta\textbf{U}_{i,j}^{n+\frac{1}{2}},\Delta\textbf{U}_{i+1,j}^{n+\frac{1}{2}},\Delta\textbf{U}_{i,j+1}^{n+\frac{1}{2}},\Delta\textbf{U}_{i+1,j+1}^{n+\frac{1}{2}})=\\	\left[\frac{S(\Delta\textbf{U}_{i,j}^{n+\frac{1}{2}})+S(\Delta\textbf{U}_{i+1,j}^{n+\frac{1}{2}})+S(\Delta\textbf{U}_{i,j+1}^{n+\frac{1}{2}})+S(\Delta\textbf{U}_{i+1,j+1}^{n+\frac{1}{2}})}{4}\right].
\end{multline*}
The forward projection step in equation \eqref{IntBL2Dc} consists of projecting the solution at time $t^{n}$ onto the staggered grid. It is performed using linear interpolations in two space dimensions in addition to Taylor expansions in space; we obtain:
\begin{multline}\label{FS2D}
	\Delta\textbf{U}_{i+\frac{1}{2},j+\frac{1}{2}}^{n} =\frac{1}{2}(\overline{\Delta\textbf{U}}_{i+\frac{1}{2},j}^n+\overline{\Delta\textbf{U}}_{i+\frac{1}{2},j+1}^n)\\-\frac{1}{16}([[\Delta\textbf{U}^{n,x}]]_{i+\frac{1}{2},j}+[[\Delta\textbf{U}^{n,x}]]_{i+\frac{1}{2},j+1}) \\
	-\frac{1}{16}([[\Delta\textbf{U}^{n,y}]]_{i,j+\frac{1}{2}}+[[\Delta\textbf{U}^{n,y}]]_{i+1,j+\frac{1}{2}}).
\end{multline}
Here, $\Delta\textbf{U}^{n,x}$ and $\Delta\textbf{U}^{n,y}$ are the spatial partial derivatives of $\Delta\textbf{U}^{n}$ that are approximated using the (MC-$\theta$) limiter \eqref{MC}.\\
Finally, the evolution step \eqref{IntBL2Dc} at time $t^{n+1}$ on the staggered nodes can be written as,
\begin{multline}\label{ES2D}
	\Delta\textbf{U}_{i+\frac{1}{2},j+\frac{1}{2}}^{n+1}=\Delta\textbf{U}_{i+\frac{1}{2},j+\frac{1}{2}}^{n}\\
	-\frac{\Delta t }{2} [D_{+}^{x}F(\Delta\textbf{U}_{i,j}^{n+\frac{1}{2}}+\tilde{\textbf{U}}_{i,j})-D_{+}^{x}F(\tilde{\textbf{U}}_{i,j})+D_{+}^{x}F(\Delta\textbf{U}_{i,j+1}^{n+\frac{1}{2}}+\tilde{\textbf{U}}_{i,j+1})\\-D_{+}^{x}F(\tilde{\textbf{U}}_{i,j+1})]\\
	-\frac{\Delta t}{2}[D_{+}^{y}G(\Delta\textbf{U}_{i,j}^{n+\frac{1}{2}}+\tilde{\textbf{U}}_{i,j})-D_{+}^{y}G(\tilde{\textbf{U}}_{i,j})+D_{+}^{y}F(\Delta\textbf{U}_{i+1,j}^{n+\frac{1}{2}}+\tilde{\textbf{U}}_{i+1,j})\\-D_{+}^{y}G(\tilde{\textbf{U}}_{i+1,j})]\\
	+\Delta t.S(\Delta\textbf{U}_{i,j}^{n+\frac{1}{2}},\Delta\textbf{U}_{i+1,j}^{n+\frac{1}{2}},\Delta\textbf{U}_{i,j+1}^{n+\frac{1}{2}},\Delta\textbf{U}_{i+1,j+1}^{n+\frac{1}{2}}).
\end{multline}
Here $D_{+}^{x}$ and $D_{+}^{y}$ are the forward differences given by,\\
$D_{+}^{x}F(\textbf{U}_{i,j})=\frac{F(\textbf{U}_{i+1,j})-F(\textbf{U}_{i,j})}{\Delta x}, D_{+}^{y}F(\textbf{U}_{i,j})=\frac{F(\textbf{U}_{i,j+1})-F(\textbf{U}_{i,j})}{\Delta y}.$\\
The predicted values in equation \eqref{ES2D} are generated at time $t^{n+\frac{1}{2}}$ using a first order Taylor expansion in time in addition to the balance law \eqref{MHDBL}:
\begin{align}\label{PS2D}
	\Delta \textbf{U}_{i,j}^{n+\frac{1}{2}}=\Delta \textbf{U}_{i,j}^{n}+\frac{\Delta t}{2}\left[-\frac{(F_{i,j}^n)'}{\Delta x}+\frac{\tilde{F}_{i,j}'}{\Delta x}-\frac{(G_{i,j}^n)'}{\Delta y}+\frac{\tilde{G}_{i,j}'}{\Delta y}+S_{i,j}^n\right],
\end{align}
where $\frac{(F_{i,j}^n)'}{\Delta x},\frac{\tilde{F}_{i,j}'}{\Delta x},\frac{(G_{i,j}^n)'}{\Delta y} $ and $\frac{\tilde{G}_{i,j}'}{\Delta y}$ denote the approximate flux derivatives with \\
$(F_{i,j}^n)' =J_{F_{i,j}^n}\cdot\textbf{U}_{i,j}^{n,x}$, $\tilde{F}_{i,j}' =J_{\tilde{F}_{i,j}}\cdot\tilde{\textbf{U}}_{i,j}^{x}$, $(G_{i,j}^n)' =J_{G_{i,j}^n}\cdot\textbf{U}_{i,j}^{n,y}$, $\tilde{G}_{i,j}' =J_{\tilde{G}_{i,j}}\cdot\tilde{\textbf{U}}_{i,j}^{y}$. Here, we also use the (MC-$\theta$) limiter \eqref{MC} to compute the slopes $\textbf{U}_{i,j}^{n,x}$, $\tilde{\textbf{U}}_{i,j}^{x}$, $\textbf{U}_{i,j}^{n,y}$, and $\tilde{\textbf{U}}_{i,j}^{y}$ in order to avoid spurious oscillations. $S_{i,j}^n$ is the discrete source term.\\
In order to retrieve the solution at the time $t^{n+1}$  on the original cells $C_{i,j}$, we project the solution obtained on the ghost cells $(\Delta\textbf{U}_{i+\frac{1}{2},j+\frac{1}{2}}^{n+1})$ back onto the oroginal grid via linear interpolations in two space dimensions and Taylor expnsions in space,
\begin{multline}\label{BS2D}
	\Delta\textbf{U}_{i,j}^{n+1}=\frac{1}{2}(\overline{\Delta\textbf{U}}_{i,j-\frac{1}{2}}^{n+1}+\overline{\Delta\textbf{U}}_{i,j+\frac{1}{2}}^{n+1})\\
	-\frac{1}{16}([[\Delta\textbf{U}^{n+1,x}]]_{(i),j-\frac{1}{2}}+[[\Delta\textbf{U}^{n+1,x}]]_{(i),j+\frac{1}{2}})\\
	-\frac{1}{16}([[\Delta\textbf{U}^{n+1,y}]]_{i-\frac{1}{2},(j)}+[[\Delta\textbf{U}^{n+1,y}]]_{i+\frac{1}{2},(j)}),
\end{multline}
where $\Delta\textbf{U}_{i,j}^{n+1,x}$ and $\Delta\textbf{U}_{i,j}^{n+1,y}$ denote the spatial partial derivatives of the numerical solution obtained at time $t^{n+1}$ and at the node $(x_i,y_j)$ approximated using the (MC-$\theta$) limiter \eqref{MC}.\\
%
%
%
To complete the presentation of the numerical scheme, we need to verify the well-balanced property of the proposed scheme and to show that it is capable of maintaining stationary solutions of the Euler system with gravitational source term.\\
Suppose that the numerical solution obtained at time $t=t^n$ satisfies $\textbf{U}_{i,j}^{n}=\tilde{\textbf{U}}_{i,j}$, i.e., $\Delta \textbf{U}_{i,j}^{n}=0.$ Performing one iteration using the proposed numerical scheme, one can show that:
\begin{enumerate}
	\item $\Delta \textbf{U}_{i,j}^{n+\frac{1}{2}}=0.$
	\item $\Delta \textbf{U}_{i+\frac{1}{2},j+\frac{1}{2}}^{n+1}=0.$
	\item $\Delta \textbf{U}_{i,j}^{n+1}=0.$\
\end{enumerate}
In fact,  it is straight forward to establish 2 and 3 once 1 is established. We will present the proof of 1 only. \\
The prediction step \eqref{PS2D} leads to 
\begin{multline}
	\Delta \textbf{U}_{i,j}^{n+\frac{1}{2}}=\Delta \textbf{U}_{i,j}^{n}+\frac{\Delta t}{2}\Bigg[ -\frac{F'(\Delta\textbf{U}_{i,j}^{n}+\tilde{\textbf{U}}_{i,j})}{\Delta x}+\frac{F'(\tilde{\textbf{U}}_{i,j})}{\Delta x}\\-\frac{G'(\Delta\textbf{U}_{i,j}^{n}+\tilde{\textbf{U}}_{i,j})}{\Delta y}+\frac{G'(\tilde{\textbf{U}}_{i,j})}{\Delta y}
	+S(\Delta\textbf{U}_{i,j}^{n},x,y)\Bigg].
\end{multline}
But since $\Delta \textbf{U}_{i,j}^{n}=0$, then we obtain, 
\begin{align*}
	\Delta \textbf{U}_{i,j}^{n+\frac{1}{2}}=\frac{\Delta t}{2}\left[-\frac{F'(\tilde{\textbf{U}}_{i,j})}{\Delta x}+\frac{F'(\tilde{\textbf{U}}_{i,j})}{\Delta x}-\frac{G'(\tilde{\textbf{U}_{i,j}})}{\Delta y}+\frac{G'(\tilde{\textbf{U}}_{i,j})}{\Delta y}\right].
\end{align*}
Hence, $\Delta \textbf{U}_{i,j}^{n+\frac{1}{2}}=0.$ Therefore, we conclude that the updated numerical solution remains stationary up to machine precision. 

\section{The constrained transport method (CTM)}
\label{sec::CTM}
In this work we consider the version of CTM developed in \cite{Touma2010}. At the end of each iteration, we apply the CTM corrections to the magnetic field components. Starting from a magnetic field that satisfies the divergence-free constraint $\nabla\cdot\textbf{B}_{i,j}^{n}=0$, we would like to prove $\nabla\cdot\textbf{B}_{i,j}^{n+1}=0$.
The discrete divergence using centered differences at time $t^n$ is given by,
\begin{align*}
	\nabla\cdot\textbf{B}_{i,j}^{n}&=\left(\frac{\partial B_x}{\partial x}\right)_{i,j}^{n}+\left(\frac{\partial B_y}{\partial y}\right)_{i,j}^{n}\\
	&=\frac{\left(B_x\right)_{i+1,j}^{n}-\left(B_x\right)_{i-1,j}^{n}}{2\Delta x}+\frac{\left(B_y\right)_{i,j+1}^{n}-\left(B_y\right)_{i,j-1}^{n}}{2\Delta y}\\
	&=0.
\end{align*}
The vector of conserved variables $\textbf{U}^{n+1}$ is computed by the numerical scheme, but $\nabla\cdot\textbf{B}_{i,j}^{n+1}$ might not be zero. Therefore, whenever needed, we correct the components of the magnetic field  $\textbf{B}_{i,j}^{n+1}$ by discretizing the induction equation at the cell centers of $C_{i,j}$,
\begin{align*}
	\frac{\partial}{\partial t}\left(
	\begin{array}{ c c}
		B_x\\
		B_y
	\end{array}\right)
	-\frac{\partial}{\partial x}\left(
	\begin{array}{ c c}
		0\\
		\Omega
	\end{array}\right)
	+\frac{\partial}{\partial y}\left(
	\begin{array}{ c c}
		\Omega\\
		0
	\end{array}\right)
	=0,
\end{align*}
where $\Omega=(-\textbf{u}\times \textbf{B})_z=-u_xB_y+u_yB_x.$ Hence, the discretization of the induction equation is the following,
$$\begin{cases}
	\frac{\left(B_x\right)_{i+\frac{1}{2},j+\frac{1}{2}}^{n+1}-\left(B_x\right)_{i+\frac{1}{2},j+\frac{1}{2}}^{n}}{\Delta t}+\frac{\Omega_{i+\frac{1}{2},j+\frac{3}{2}}^{n+\frac{1}{2}}-\Omega_{i+\frac{1}{2},j-\frac{1}{2}}^{n+\frac{1}{2}}}{2\Delta y}=0,\\
	\frac{\left(B_y\right)_{i+\frac{1}{2},j+\frac{1}{2}}^{n+1}-\left(B_y\right)_{i+\frac{1}{2},j+\frac{1}{2}}^{n}}{\Delta t}-\frac{\Omega_{i+\frac{3}{2},j+\frac{1}{2}}^{n+\frac{1}{2}}-\Omega_{i-\frac{1}{2},j+\frac{1}{2}}^{n+\frac{1}{2}}}{2\Delta x}=0.
\end{cases}$$
Then,
\begin{equation}\label{cpt}
	\begin{cases}
		\left(B_x\right)_{i+\frac{1}{2},j+\frac{1}{2}}^{n+1}=\left(B_x\right)_{i+\frac{1}{2},j+\frac{1}{2}}^{n}-\frac{\Delta t}{2\Delta y} \left(\Omega_{i+\frac{1}{2},j+\frac{3}{2}}^{n+\frac{1}{2}}-\Omega_{i+\frac{1}{2},j-\frac{1}{2}}^{n+\frac{1}{2}}\right),\\
		\left(B_y\right)_{i+\frac{1}{2},j+\frac{1}{2}}^{n+1}=\left(B_y\right)_{i+\frac{1}{2},j+\frac{1}{2}}^{n}+\frac{\Delta t}{2\Delta x}\left(\Omega_{i+\frac{3}{2},j+\frac{1}{2}}^{n+\frac{1}{2}}-\Omega_{i-\frac{1}{2},j+\frac{1}{2}}^{n+\frac{1}{2}}\right).
	\end{cases}
\end{equation}
Now, we compute $\Omega_{i+\frac{1}{2},j+\frac{1}{2}}^{n+\frac{1}{2}}$
using the numerical solution computed at time $t^n$ and $t^{n+1}$ in order to obtain second order of accuracy in time, 
\begin{align*}
	\Omega_{i+\frac{1}{2},j+\frac{1}{2}}^{n+\frac{1}{2}}&=\frac{1}{2} \left[\Omega_{i+\frac{1}{2},j+\frac{1}{2}}^{n+1}+\Omega_{i+\frac{1}{2},j+\frac{1}{2}}^{n}\right],\\
	&=\frac{1}{2}\left[\Omega_{i+\frac{1}{2},j+\frac{1}{2}}^{n+1}+\frac{\Omega_{i,j}^{n}+\Omega_{i+1,j}^{n}+\Omega_{i,j+1}^{n}+\Omega_{i+1,j+1}^{n}}{4}\right].
\end{align*}
Next, we calculate $\nabla\cdot(\textbf{B})_{i+\frac{1}{2},j+\frac{1}{2}}^{n+1}$
\begin{multline}\label{divstag}
	\nabla\cdot(\textbf{B})_{i+\frac{1}{2},j+\frac{1}{2}}^{n+1}=\frac{\left(B_x\right)_{i+\frac{3}{2},j+\frac{1}{2}}^{n+1}-\left(B_x\right)_{i-\frac{1}{2},j+\frac{1}{2}}^{n+1}}{2\Delta x}+\frac{\left(B_y\right)_{i+\frac{1}{2},j+\frac{3}{2}}^{n+1}-\left(B_y\right)_{i+\frac{1}{2},j-\frac{1}{2}}^{n+1}}{2\Delta y}.
\end{multline}
Substituting the magnetic field components on the staggered grid in \eqref{divstag} from their values in \eqref{cpt} leads to, 
\begin{multline}
	\nabla\cdot(\textbf{B})_{i+\frac{1}{2},j+\frac{1}{2}}^{n+1}=\frac{1}{4}\left[\nabla\cdot\textbf{B}_{i,j}^n+\nabla\cdot\textbf{B}_{i+1,j+1}^n+\nabla\cdot\textbf{B}_{i+1,j}^n+\nabla\cdot\textbf{B}_{i,j+1}^n\right]=0.
\end{multline}
Finally, we compute the magnetic field on the main grid $\textbf{B}_{i,j}^{n+1}$ as the average of its values on the staggered grid,
\begin{align*}
	\textbf{B}_{i,j}^{n+1}=\frac{1}{4}\left[\textbf{B}_{i+\frac{1}{2},j+\frac{1}{2}}^{n+1}+\textbf{B}_{i+\frac{1}{2},j-\frac{1}{2}}^{n+1}+\textbf{B}_{i-\frac{1}{2},j+\frac{1}{2}}^{n+1}+\textbf{B}_{i-\frac{1}{2},j-\frac{1}{2}}^{n+1}\right].
\end{align*}
Hence,
\begin{align}
	\nabla\cdot\textbf{B}_{i,j}^{n+1}=0.
\end{align}

\section{Numerical Experiments}
\label{sec::numerical results}
A list of numerical experiments has been considered in order to verify the robustness and accuracy of our method. The time-step is computed with CFL number 0.485. The MC-$\theta$ limiter is used with $\theta=1.5$.

\subsection{2D shock tube problem}
\begin{figure}[ht]
	\centering
	\includegraphics[width=1\linewidth]{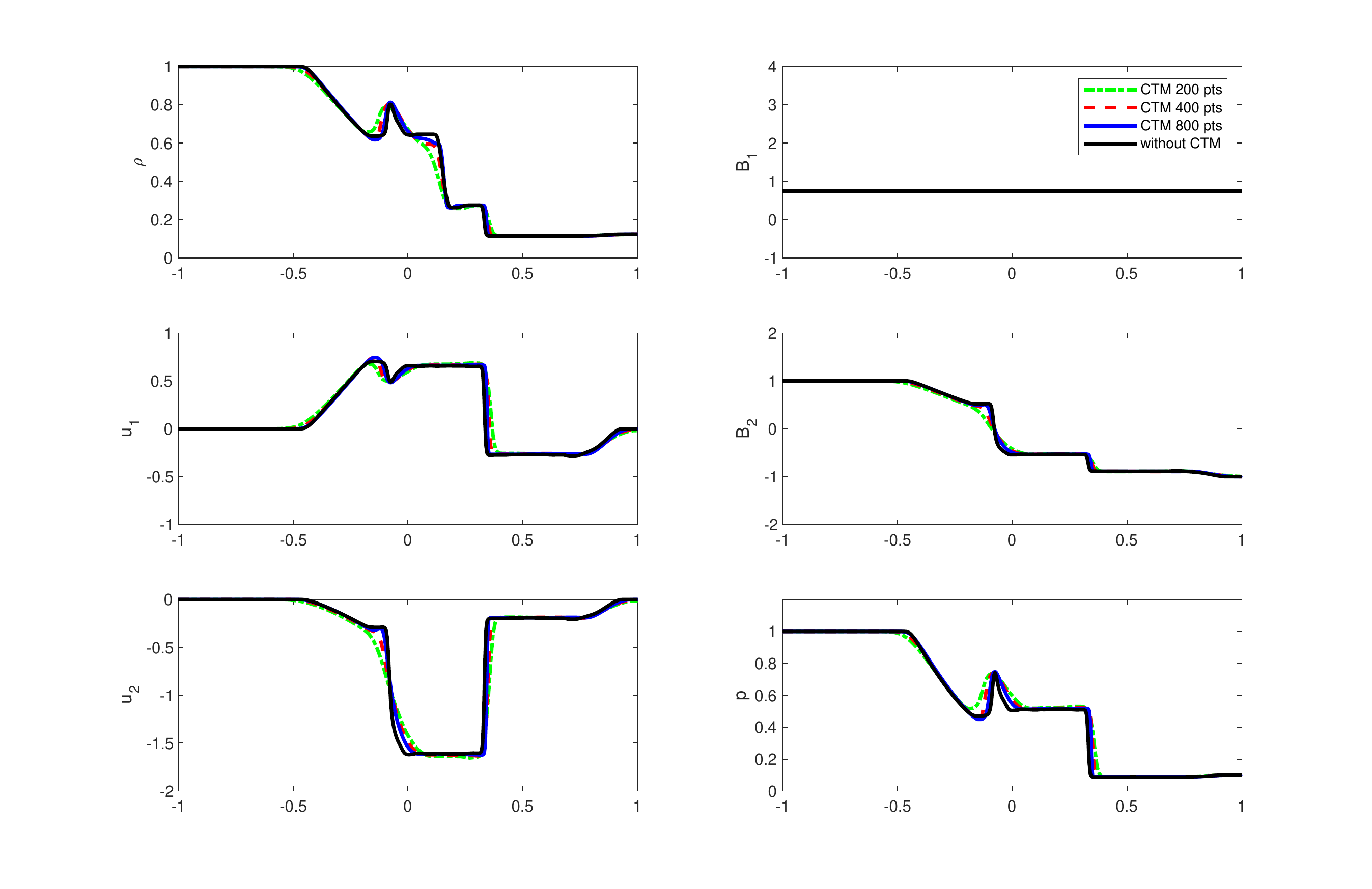}
	\caption{2D shock tube problem: a cross section of the components at time $t=0.25$ with applying CTM on different grids and without applying CTM.}
	\label{RP2Dcomparison}
\end{figure}

For the first numerical test case, we consider the Brio–Wu shock tube problem for the system of ideal MHD equations extracted from \cite{ArminjonTouma2005}. The simulation takes place over the computational domain $[-1,1]\times[-1,1]$. 
$U=[\rho,u_1,u_2,u_3,B_{2},B_{3},p]$ is initially given  as $U=[1,0,0,0,1,0,1]$ for $x<0$ and $U=[0.125,0,0,0,-1,0,0.1]$ for $x>0$ and $B_{1}=0.75$. This test case features seven discontinuities. We compute the solution at the final time $t=0.25$ with and without applying CTM in figure \ref{RP2Dcomparison}. The cross sections show a very good agreement with the results in the literature. In order to investigate the effect of the CTM on the computed solution, we did a convergence study in figure \ref{RP2Dcomparison} while applying the CTM. As it is very clear in the figure above, applying the CTM for the UC schemes has a small smearing out effect on the solution.

\subsection{Four stages Ideal MHD Riemann problem}
This test case is considered to prove the ability of our scheme to solve ideal MHD problems and preserve the divergence-free constraint. The initial data consist of four constant states \cite{ArminjonTouma2005, Touma2010}. The initial four constant states are given as follows, 
\begin{equation}
	(\rho,u_1,u_2,p)=
	\begin{cases}
		(1,0.75,0.5,1) \hspace{0.5 cm} \text{if}\hspace{0.1 cm}  x>0 \hspace{0.1 cm}  \text{and}\hspace{0.1 cm}  y>0\\
		(2,0.75,0.5,1) \hspace{0.5 cm} \text{if}\hspace{0.1 cm}  x<0 \hspace{0.1 cm}  \text{and}\hspace{0.1 cm}  y>0\\
		(1,-0.75,0.5,1) \hspace{0.5 cm} \text{if}\hspace{0.1 cm}  x<0 \hspace{0.1 cm}  \text{and}\hspace{0.1 cm}  y<0\\
		(3,-0.75,-0.5,1) \hspace{0.5 cm} \text{if}\hspace{0.1 cm}  x>0 \hspace{0.1 cm}  \text{and}\hspace{0.1 cm}  y<0\\
	\end{cases}
\end{equation}
with an initial uniform magnetic field $\textbf{B}=(2,0,1)$. The numerical solution is computed in the square $\left[-1,1\right]\times\left[-1,1\right]$ on 400$\times$400 grid points.
\begin{figure}[ht]
	\centering
	\begin{subfigure}[b]{0.5\linewidth}
		\centering
		\includegraphics[width=1\linewidth]{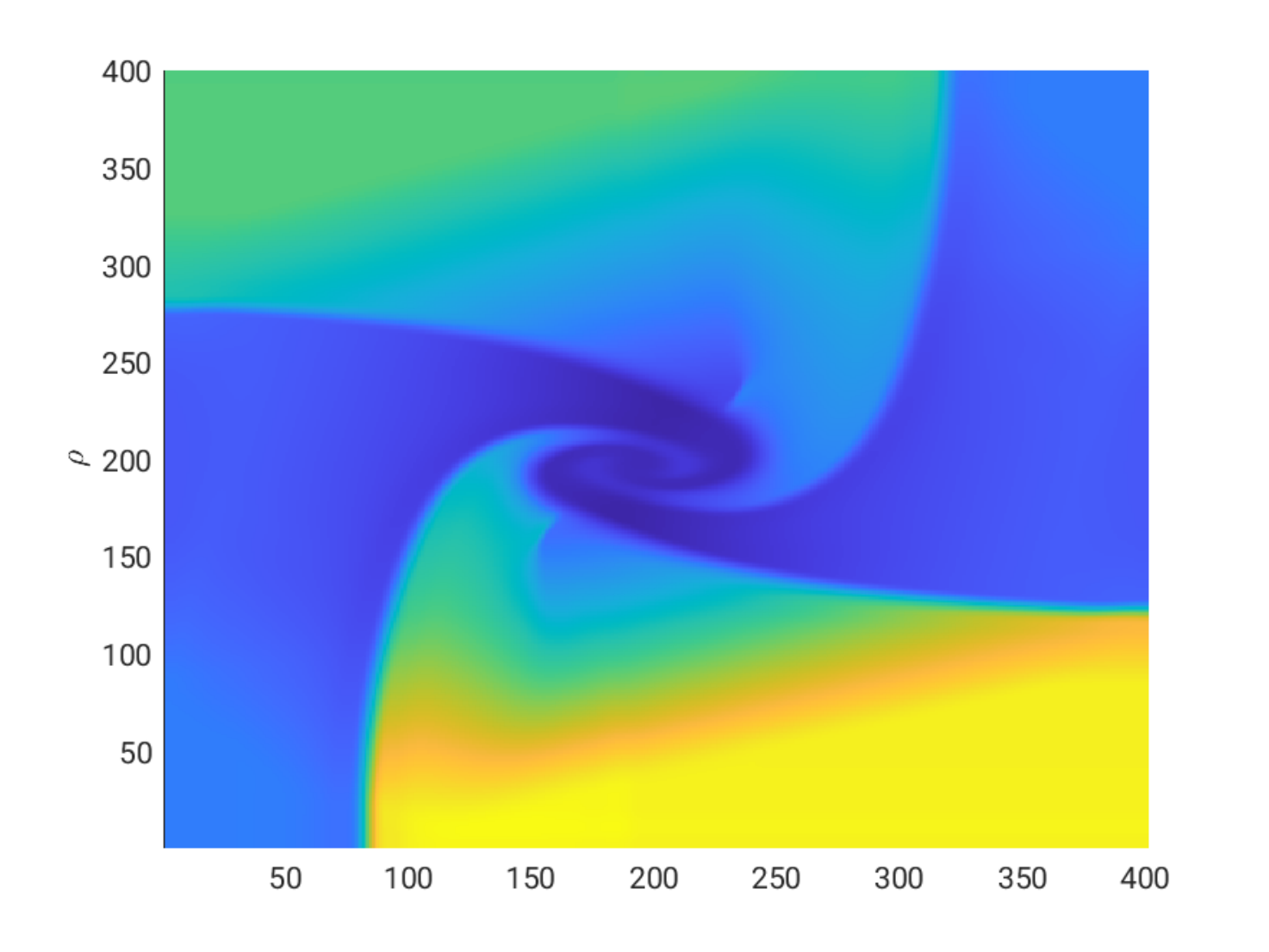}
	\end{subfigure}%
	\begin{subfigure}[b]{0.5\linewidth}
		\centering
		\includegraphics[width=1\linewidth]{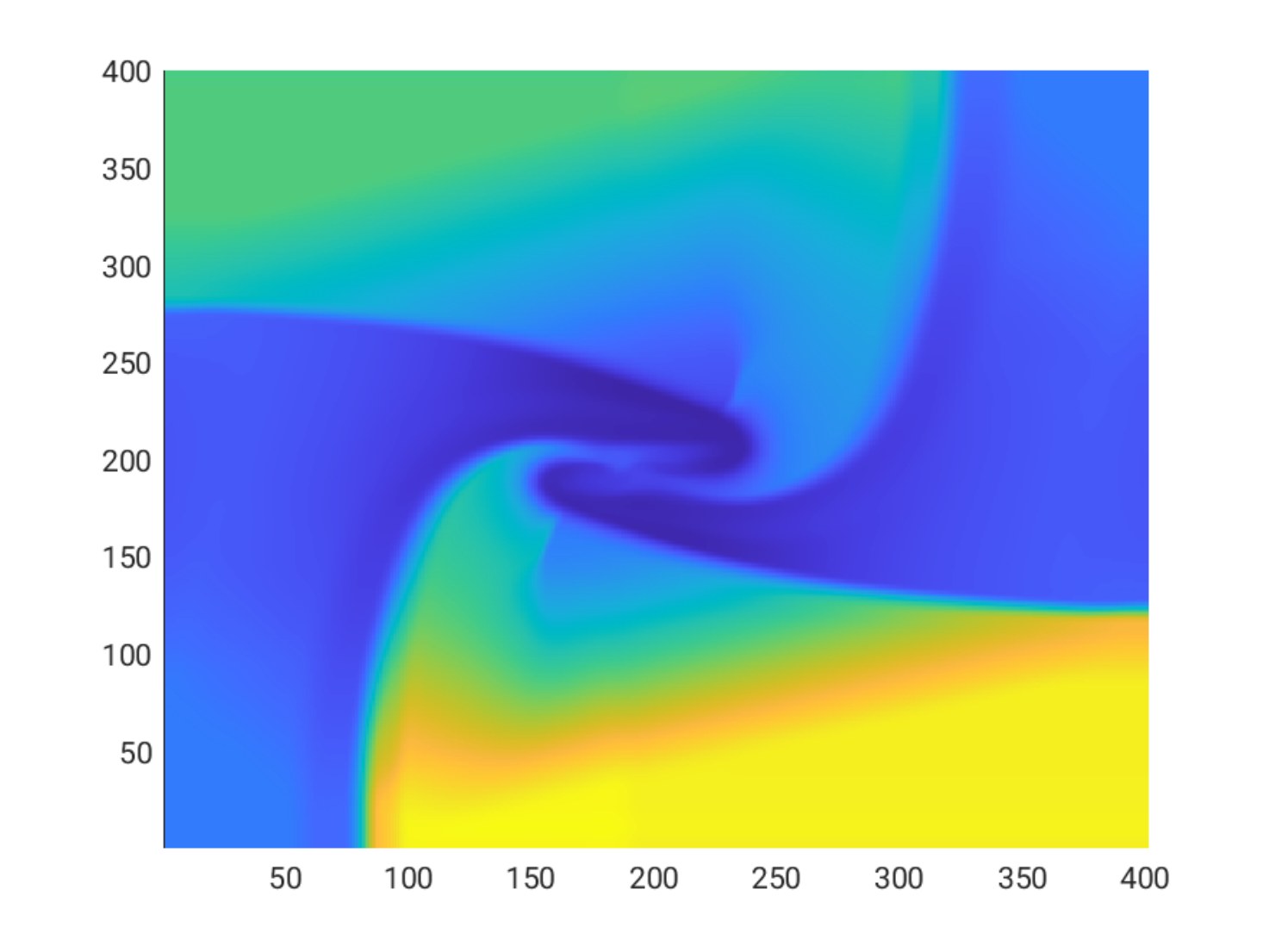}
	\end{subfigure}
	\caption{Four stages Riemann problem: $\rho$ with CTM (left) and without CTM (right) at the final time $t=0.8$.}
	\label{4rprho}
\end{figure}

\begin{figure}[ht]
	\centering
	\begin{subfigure}[b]{0.5\linewidth}
		\centering
		\includegraphics[width=1\linewidth]{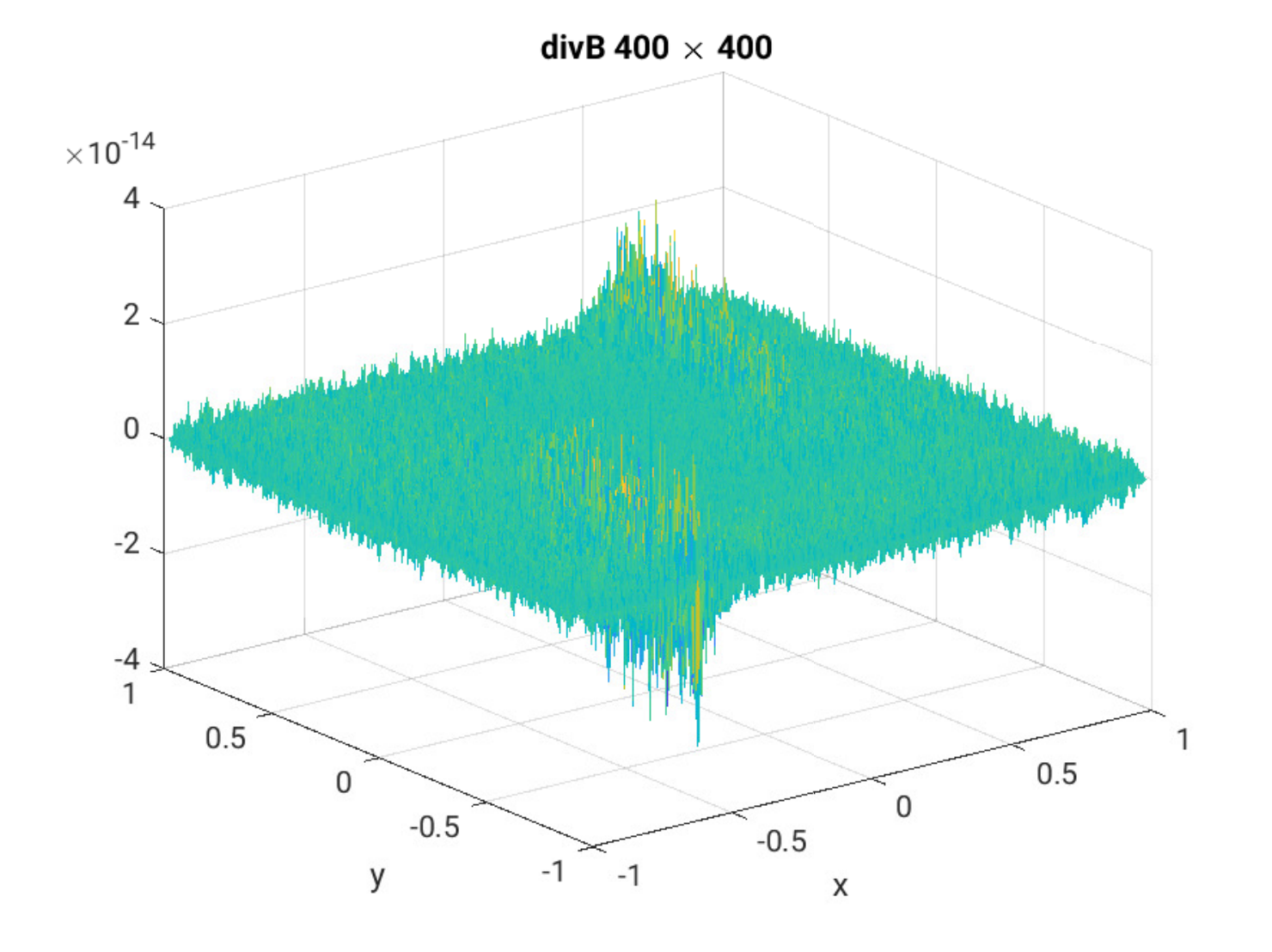}
	\end{subfigure}%
	\begin{subfigure}[b]{0.5\linewidth}
		\centering
		\includegraphics[width=1\linewidth]{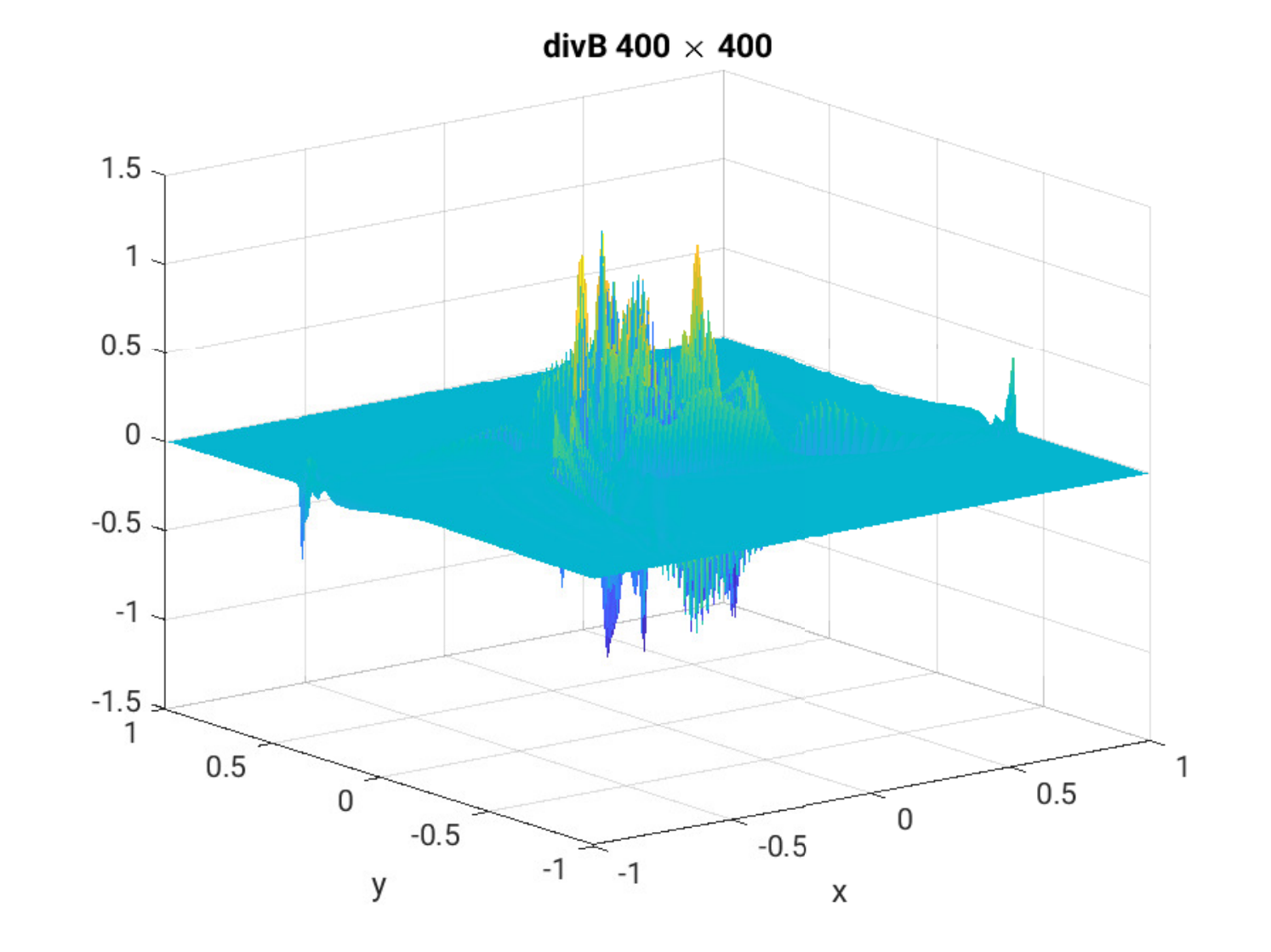}
	\end{subfigure}
	\caption{Four stages Riemann Problem: div$\textbf{B}$ with CTM (left) and without CTM (right) at the final time $t=0.8$.}
	\label{4rpdiv}
\end{figure}
Figure \ref{4rprho} illustrates the density at the final time $t_f=0.8$ with and without applying constrained transport treatment to the magnetic field components. Similar comparison on the divergence of the magnetic field is illustrated in figure \ref{4rpdiv}. The results highlight the robustness of the numerical scheme in the sense that even without treatment we are able to show numerical simulation while other schemes simply blow up without special treatment of the magnetic field.

\subsection{MHD vortex}
\begin{figure}[ht]
	\centering
	\begin{subfigure}[b]{0.3\linewidth}
		\centering
		\includegraphics[width=1\linewidth]{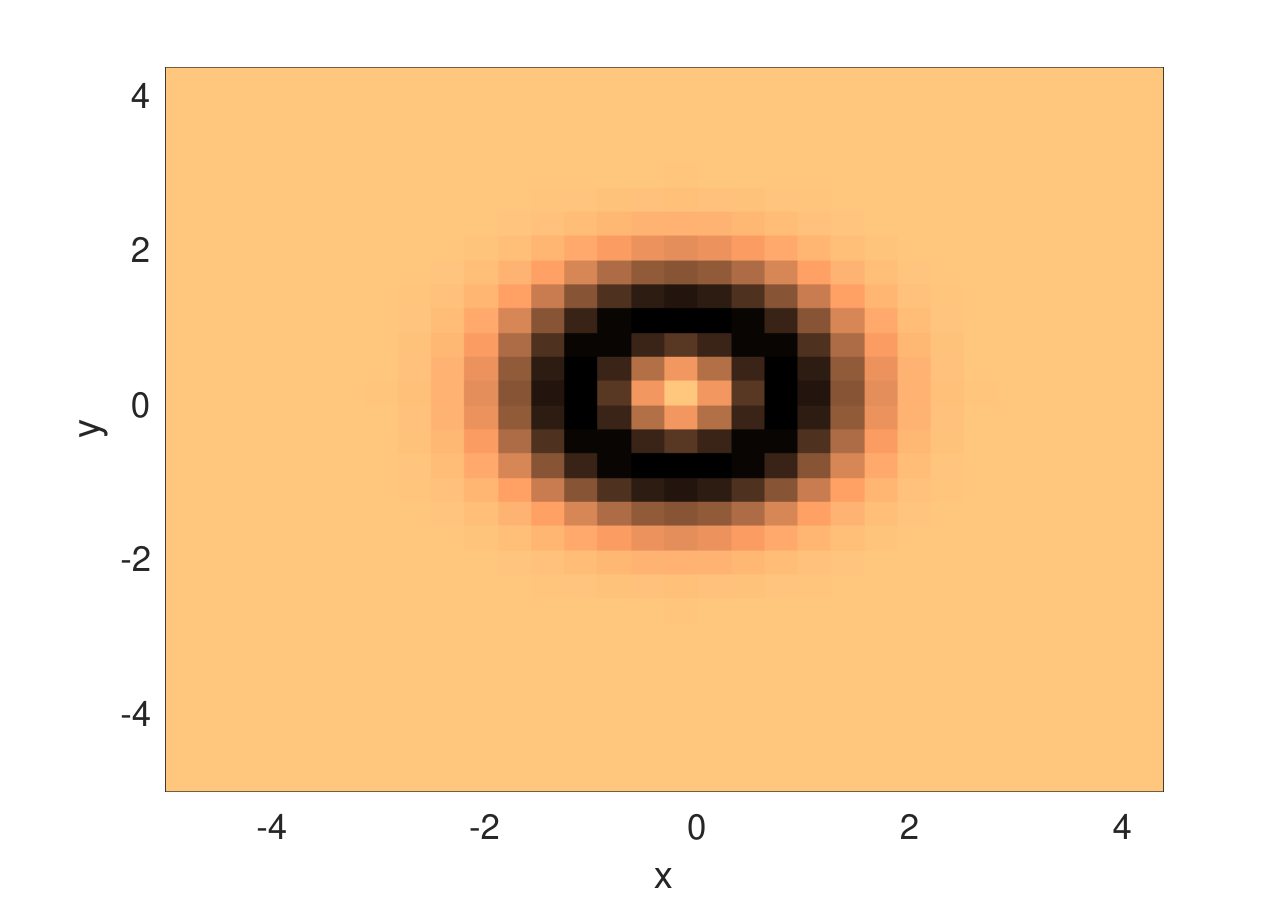}
		\subcaption*{$32 \times 32$}
	\end{subfigure}%
	\begin{subfigure}[b]{0.3\linewidth}
		\centering
		\includegraphics[width=1\linewidth]{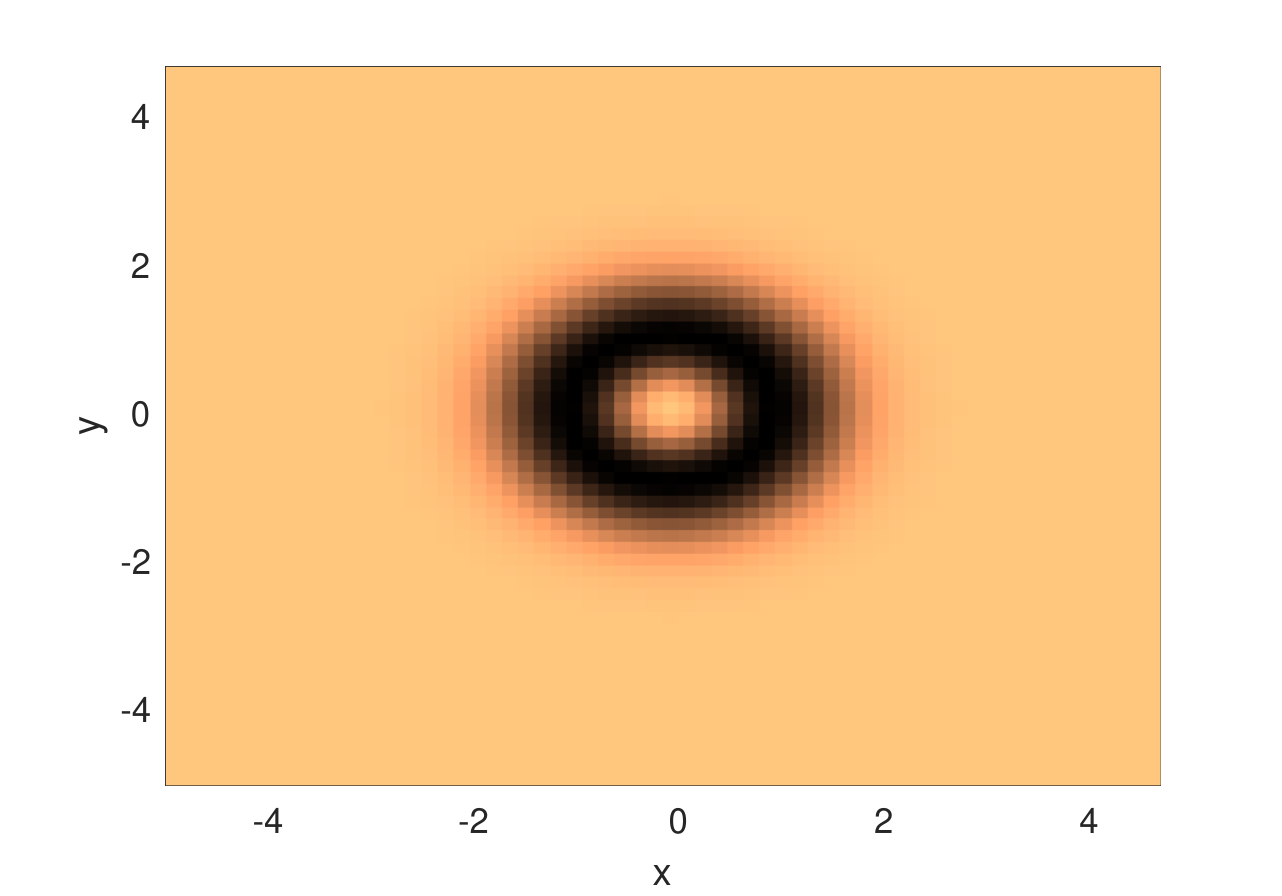}
		\subcaption*{$64 \times 64$}
	\end{subfigure}%
	\begin{subfigure}[b]{0.3\linewidth}
		\centering
		\includegraphics[width=1\linewidth]{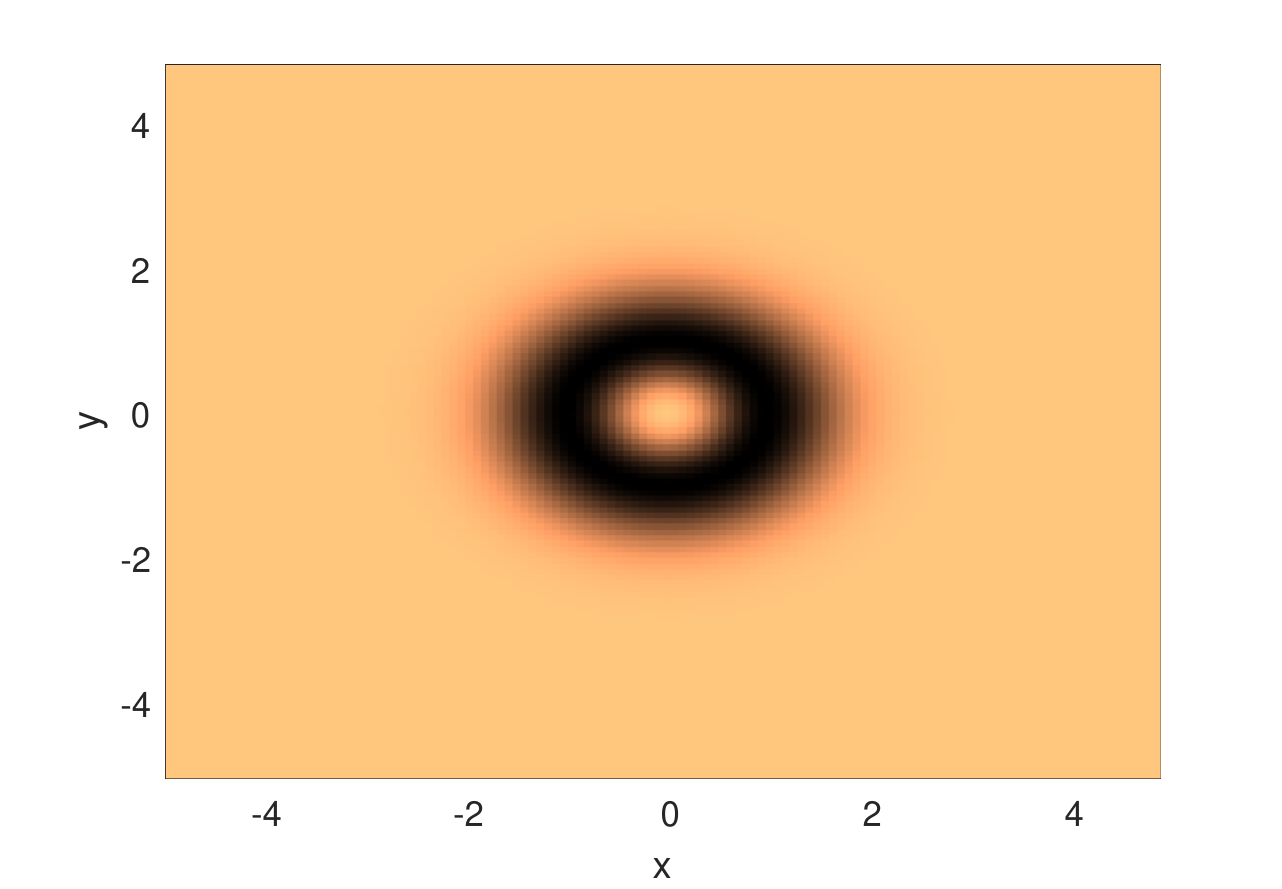}
		\subcaption*{$128 \times 128$}
	\end{subfigure}
	\caption{MHD vortex: pressure profile at the final time on different grid points.}
	\label{mhdvortex}
\end{figure}

For our third test case, we consider the MHD vortex for the homogeneous ideal MHD equations \cite{jonas}. The initial data represent a moving stationary solution of the system of the ideal MHD equations and are given by,
$r^2=x^2+y^2$,
$\rho=1$,
$u_1=u_0-\kappa_p\exp(\frac{1-r^2}{2})y$,
$u_2=v_0+\kappa_p\exp(\frac{1-r^2}{2})x$,
$u_3=0$,
$B_1=-m_p\exp(\frac{1-r^2}{2})y$,
$B_2=-m_p\exp(\frac{1-r^2}{2})x$, 
$B_3=0$, and
$p=1+\left(\frac{m_p^2}{2}(1-r^2)-\frac{\kappa_p^2}{2}\right)$.
We set the parameters $m_p=1, \kappa_p=1, u_0=0$, and $v_0=0$. The vortex is advected through the domain $\left[-5,5\right]\times\left[-5,5\right]$ with a velocity $(u_0,v_0)$. Steady state boundary conditions are used in this test case.
In figure \ref{mhdvortex}, we present the pressure profile at the final time $t=100\frac{2\pi}{\sqrt{e}\kappa_p}\approx100\frac{3.14}{\kappa_p}$ on different grids. The steady state gets preserved exactly as the background solution $\tilde{\textbf{U}}$ is the vortex itself. 

\subsection{Hydrodynamic wave propagation}
\begin{figure}[ht]
	\centering
	\begin{subfigure}[b]{0.5\linewidth}
		\centering
		\includegraphics[width=1\linewidth]{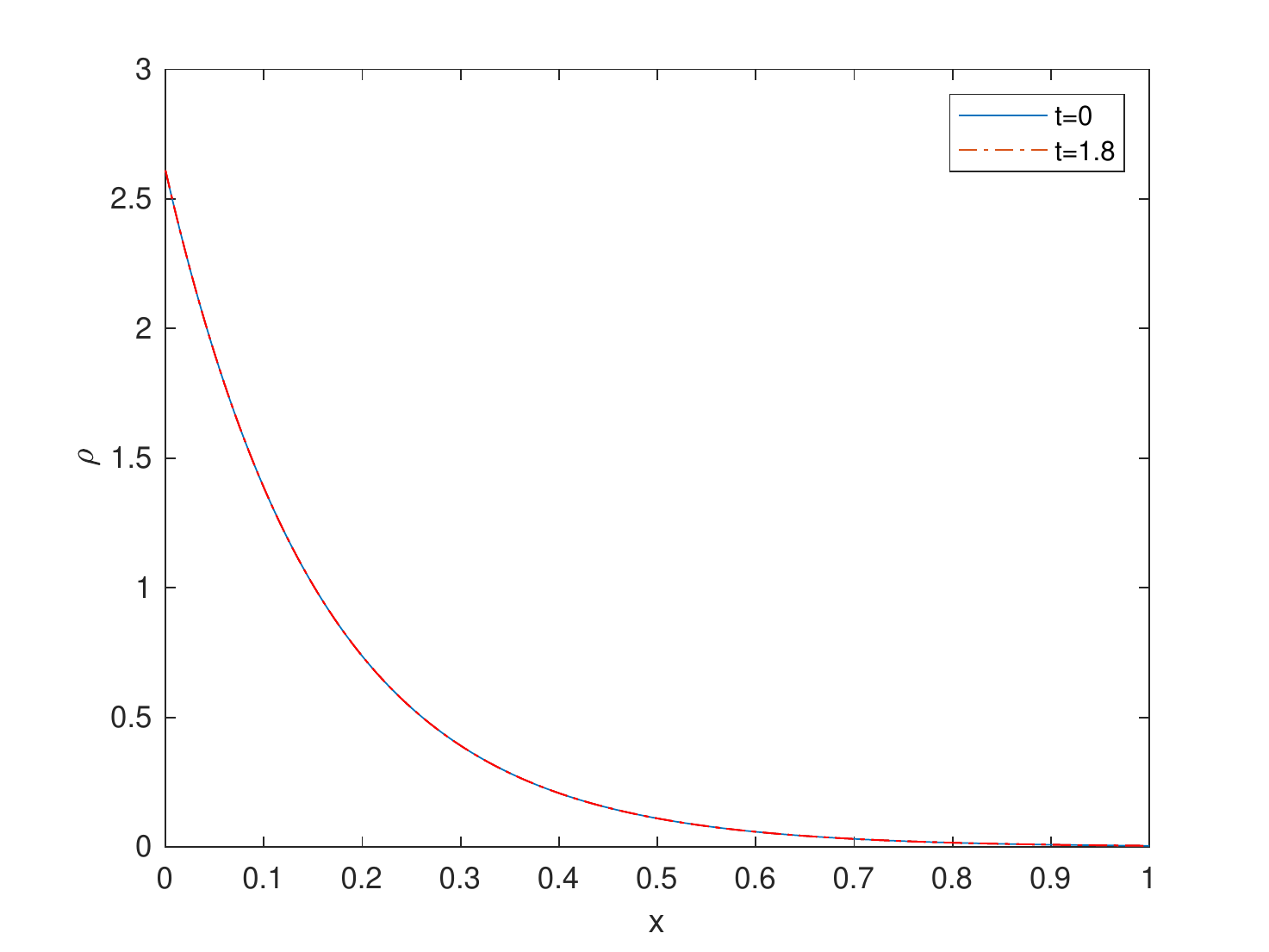}
	\end{subfigure}%
	\begin{subfigure}[b]{0.5\linewidth}
		\centering
		\includegraphics[width=1\linewidth]{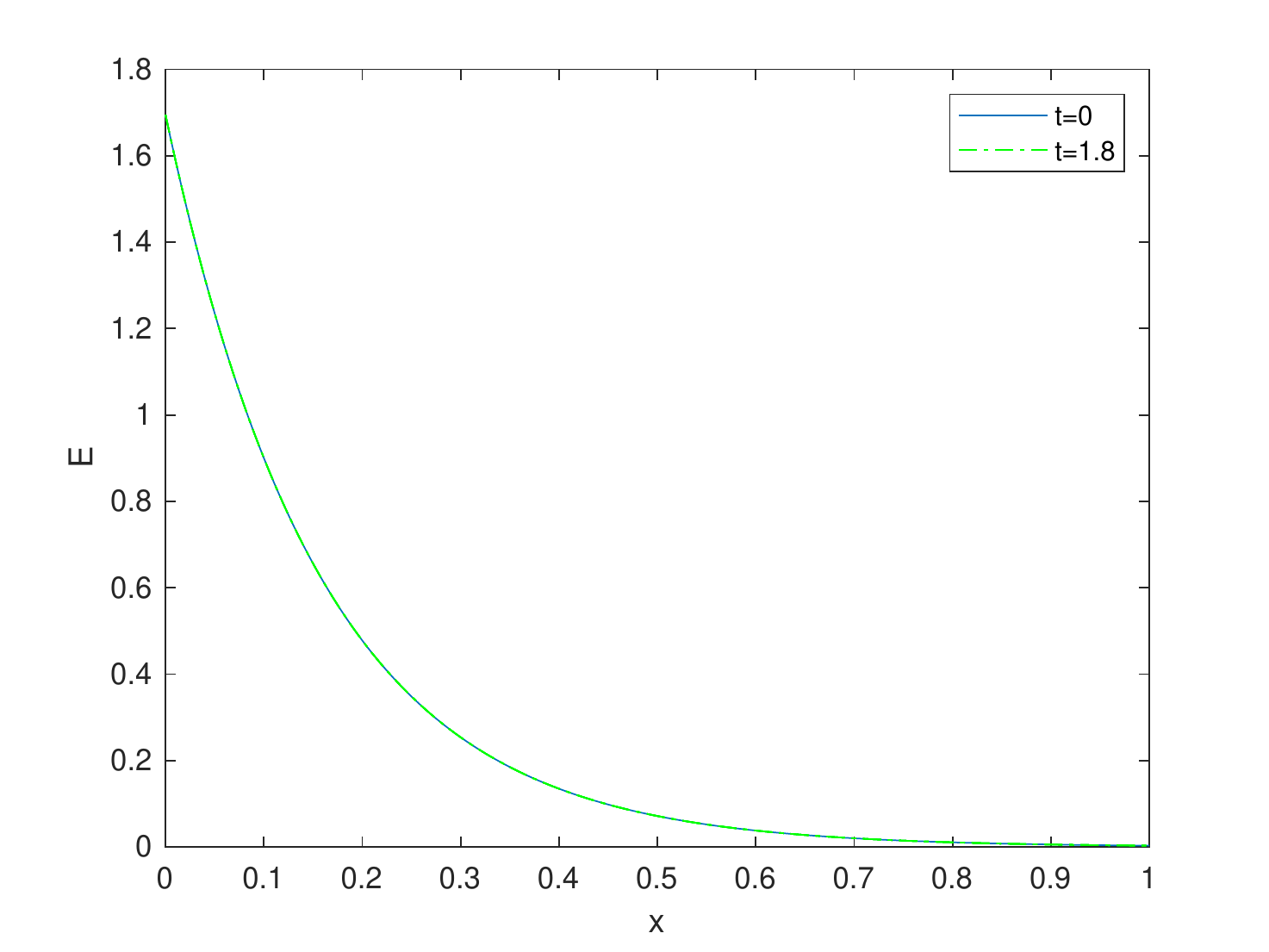}
	\end{subfigure}%
	\caption{Hydrodynamic wave propagation: a comparison of the cross sections of the density $\rho$ (left) and the energy $E$ (right) initially and at the final time $t=1.8$.}
	\label{hydrocut}
\end{figure}

\begin{figure}[ht]
	\centering
	\begin{subfigure}[b]{0.5\linewidth}
		\centering
		\includegraphics[width=1\linewidth]{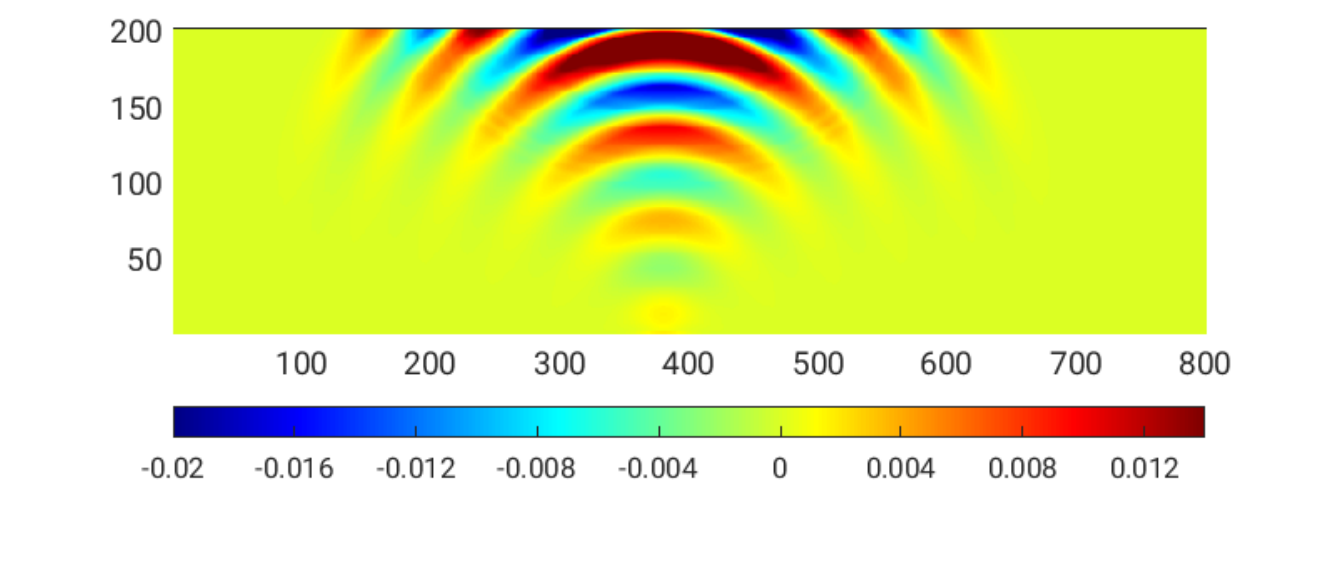}
	\end{subfigure}%
	\begin{subfigure}[b]{0.5\linewidth}
		\centering
		\includegraphics[width=1\linewidth]{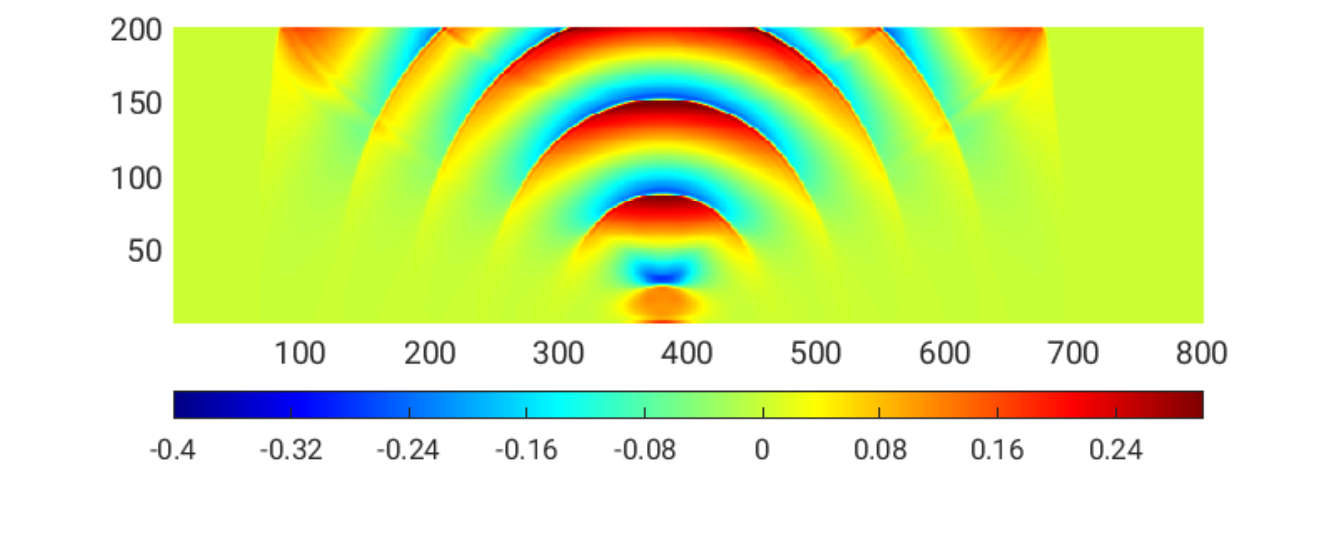}
	\end{subfigure}%
	\caption{Hydrodynamic wave propagation: wave profile $u_2$ for $c=0.003$ (left) and $c=0.3$ (right) at the final time $t=1.8$.}
	\label{hydro}
\end{figure}

The aim of this test case is to test the well-balanced property of the subtraction method by simulating a steady state solution under hydrodynamic wave propagation. The experiment is carried out in two steps. The first step is to check that the subtraction method preserves the steady state. The initial data are the hydrodynamic steady state in the computational domain $\left[0,4\right]\times\left[0,1\right]$.
\begin{align}
	\rho(x,y)=\rho_{0}\exp(-\frac{y}{H}), p(x,y)=p_{0}\exp(-\frac{y}{H}),  \textbf{u}=0, \textbf{B}=0.
\end{align}
With $H=\frac{p_0}{g\rho_{0}}=0.158$, $p_{0}=1.13$ and $g=2.74$. The subtraction method  preserves the hydrodynamic steady state exactly after choosing the reference solution  $\tilde{\textbf{U}}$ at the steady state itself. Figure \ref{hydrocut} shows a very simple comparison of the density and the energy cross sections at $t=0$ and the final time $t=1.8$. The second step is to add perturbation to the steady state as a time dependent sinusoidal wave that propagates from the bottom boundary of the vertical velocity and exits from the top one.  The wave formula is the following,
\begin{equation}
	u_{2_{i,\{0,-1\}}}^{n}=\exp(-100(x_{i,\{0,-1\}}-1.9)^2)c\sin(6\pi t^n).
\end{equation}
The bottom boundary is a localized piston at $x=1.9$. Figure \ref{hydro} shows the profile of the wave at the final time $t=1.8$ for $c=0.003$ (left) and for $c=0.3$ (right) for $800 \times 200$ grid points.  The waves propagate in both cases from bottom to top under the effect of the pressure and gravity forces. The case where $c=0.003$  models a small perturbation and  $c=0.3$ models a stronger wave. The results are in a very good agreement with the ones in \cite{Fuchs2010}. Additionally, they match the results of the most accurate (third order) of the three schemes compared in \cite{Fuchs2010}. Hence, the scheme is well-balanced in the sense that it preserves the steady state and can capture its perturbations.  

\subsection{MHD wave propagation}
\begin{figure}[ht]
	\centering
	\includegraphics[width=0.5\linewidth]{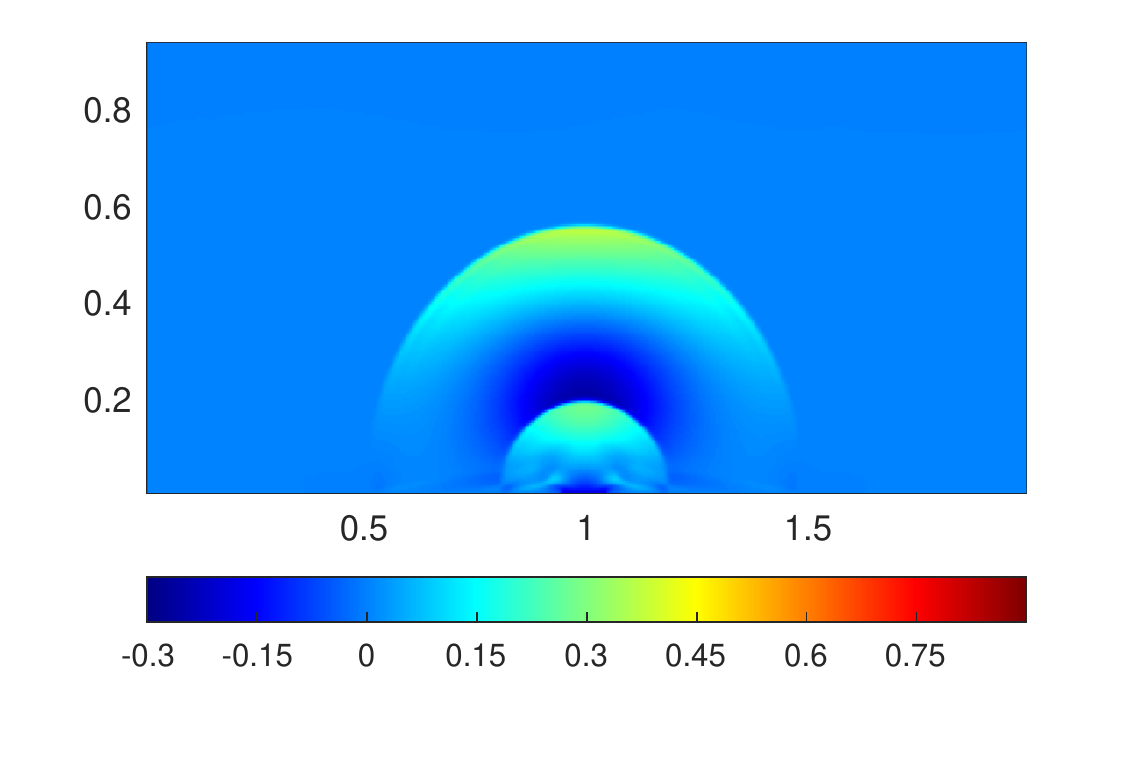}
	\caption{MHD wave propagation: velocity in a direction parallel to the magnetic field $u_B=<\textbf{u},\textbf{B}>/|\textbf{B}|$ for $\mu=0$ on 400 $\times$ 200 grid points at the final time $t=0.54$.}
	\label{magnetoB=0}
\end{figure}

\begin{figure}[ht]
	\begin{subfigure}[b]{0.45\linewidth}
		\centering
		\includegraphics[width=1\linewidth]{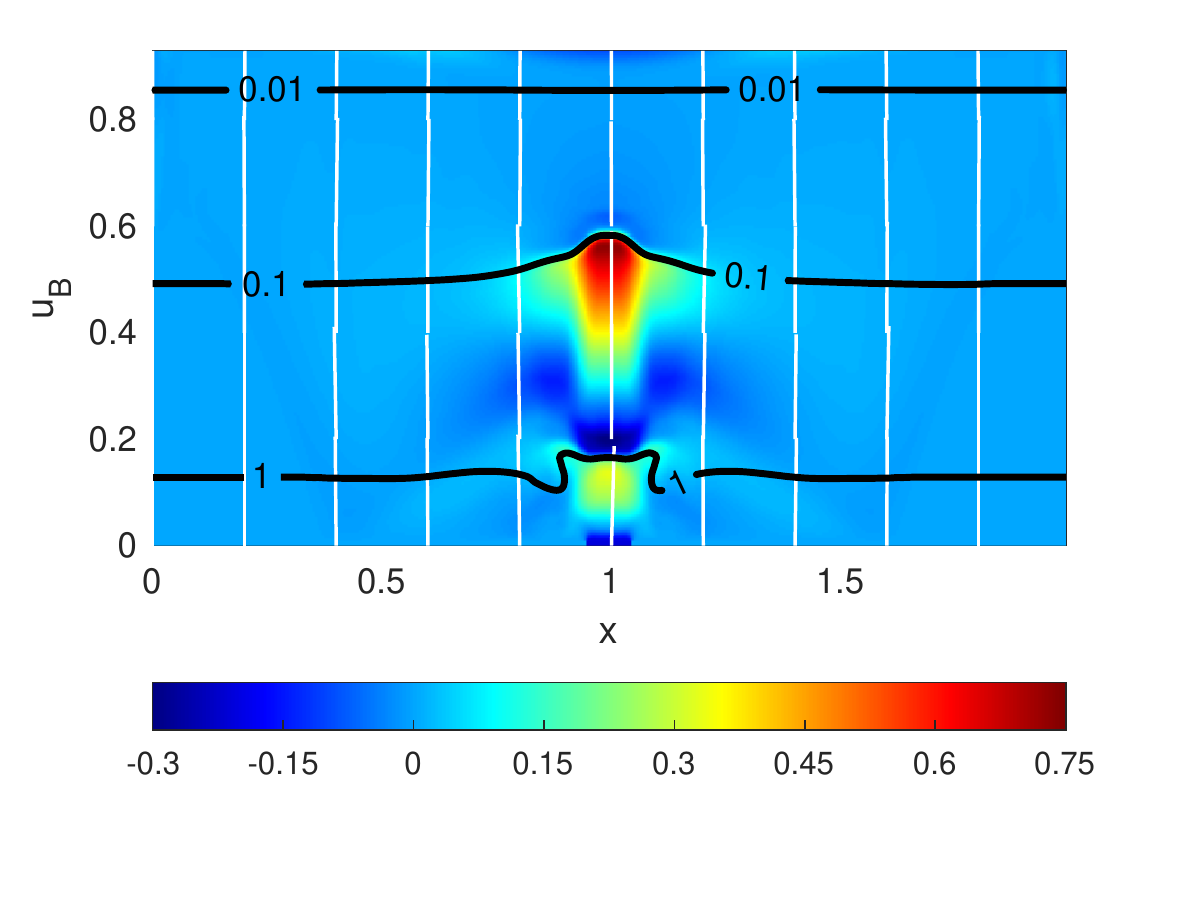}
	\end{subfigure}%
	\begin{subfigure}[b]{0.5\linewidth}
		\centering
		\includegraphics[width=1\linewidth]{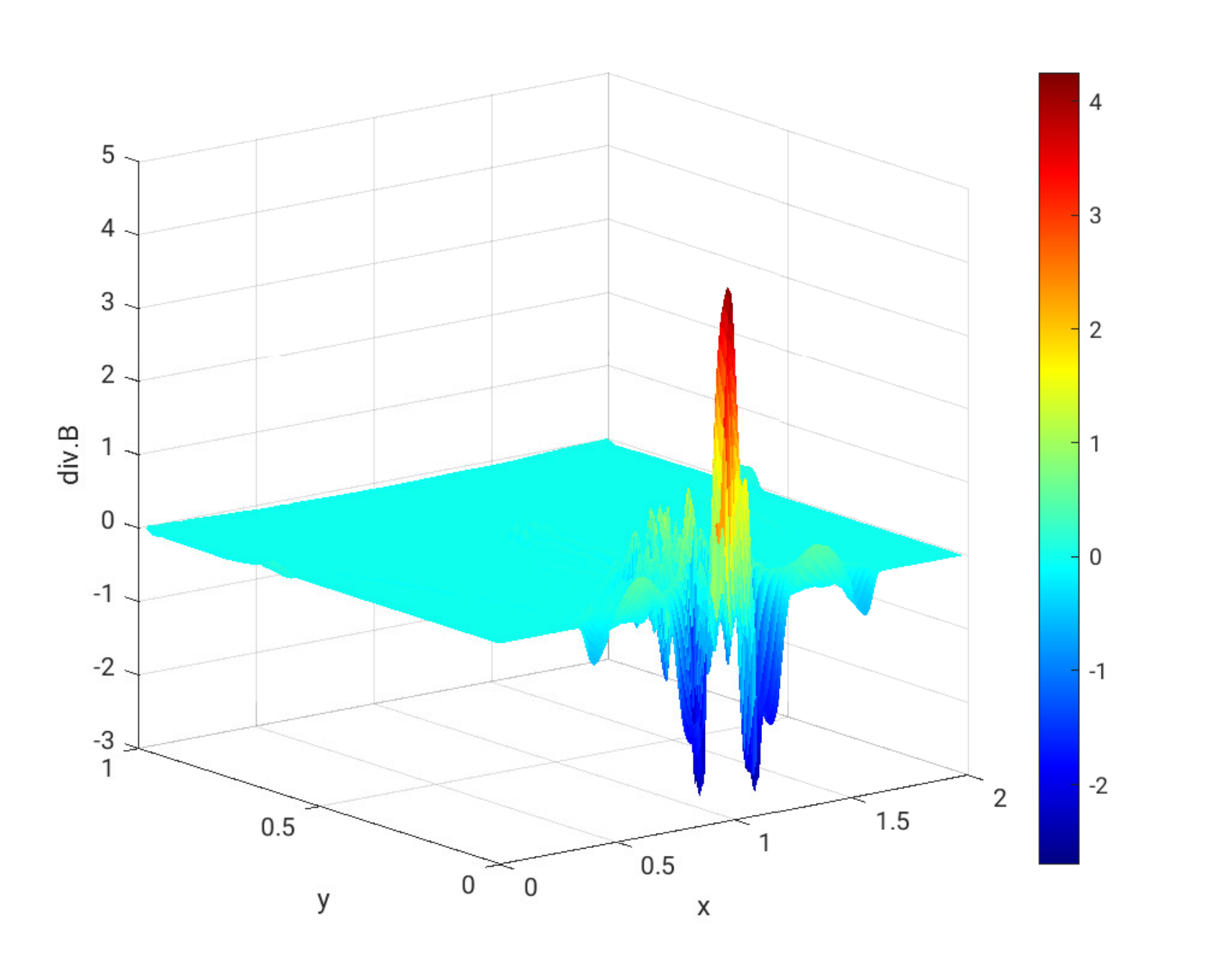}
	\end{subfigure}
	\caption{MHD wave propagation: velocity in a direction parallel to the magnetic field $u_B=<\textbf{u},\textbf{B}>/|\textbf{B}|$ for $\mu=1$ on 400 $\times$ 200 grid points at the final time $t=0.54$ without CTM.}
	\label{magnetoB=1}
\end{figure}

\begin{figure}[ht]
	\begin{subfigure}[b]{0.45\linewidth}
		\centering
		\includegraphics[width=1\linewidth]{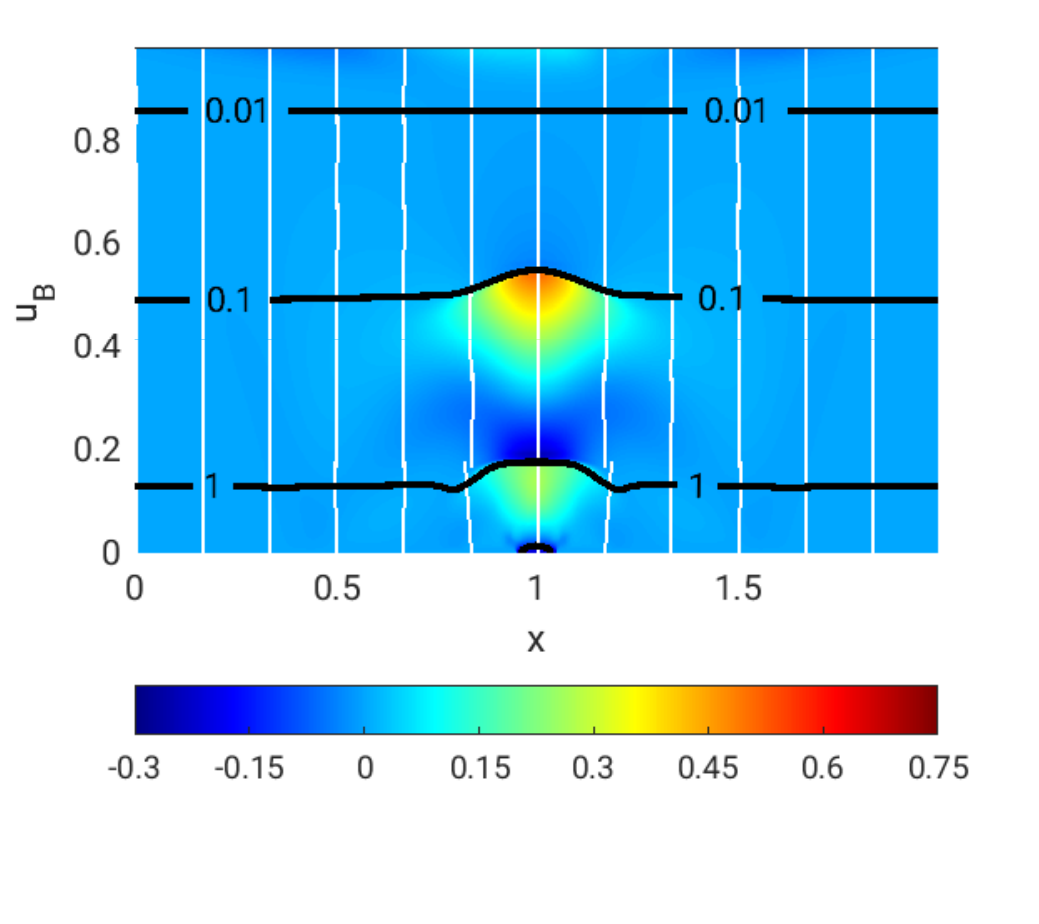}
	\end{subfigure}%
	\begin{subfigure}[b]{0.5\linewidth}
		\centering
		\includegraphics[width=1\linewidth]{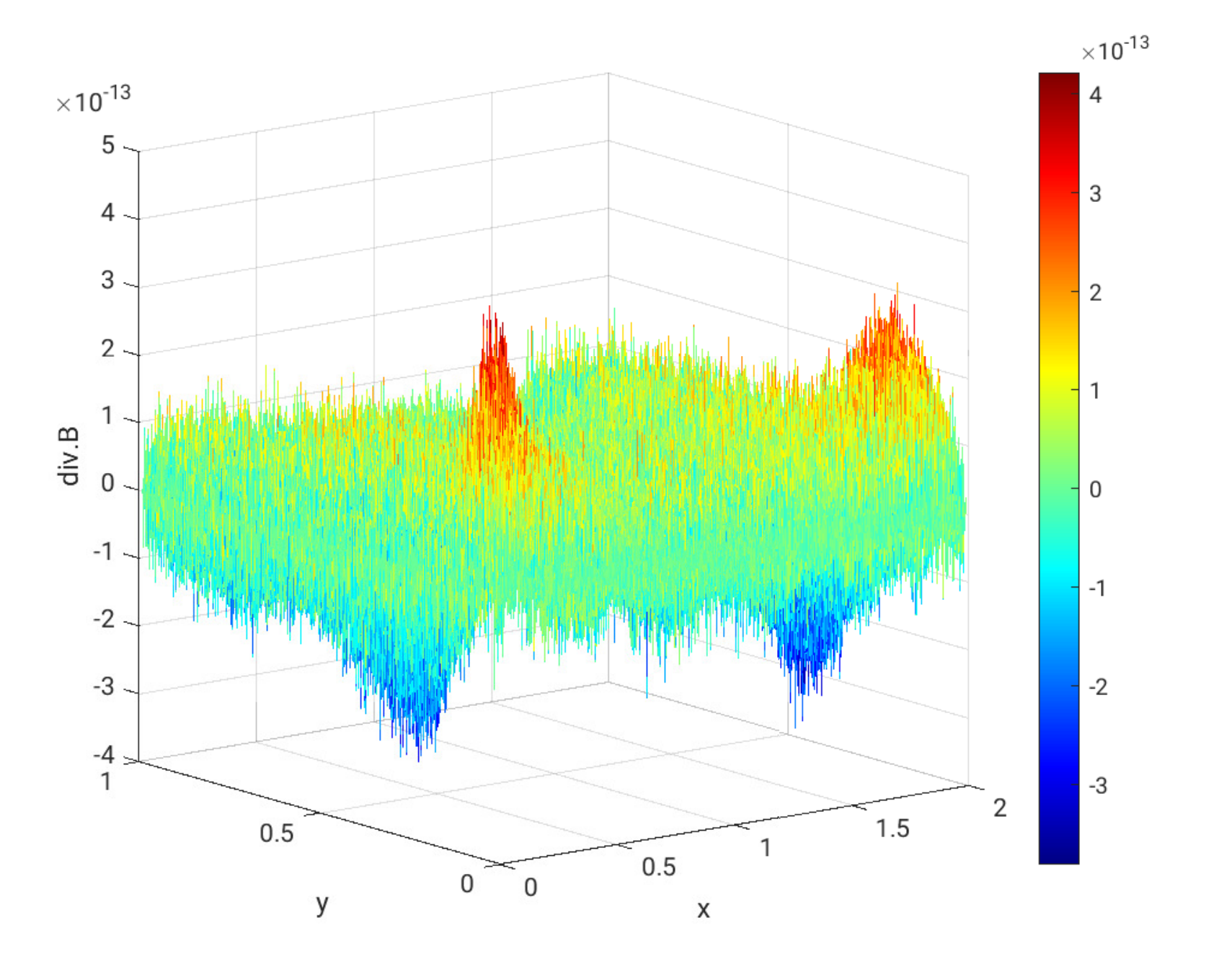}
	\end{subfigure}
	\caption{MHD wave propagation: velocity in a direction parallel to the magnetic field $u_B=<\textbf{u},\textbf{B}>/|\textbf{B}|$ for $\mu=1$ on 1200 $\times$ 600 grid points at the final time $t=0.54$ with CTM.}
	\label{magnetoB=1CTM}
\end{figure}

\begin{figure}[ht]
	\centering
	\begin{subfigure}[b]{0.3\linewidth}
		\centering
		\includegraphics[width=1\linewidth]{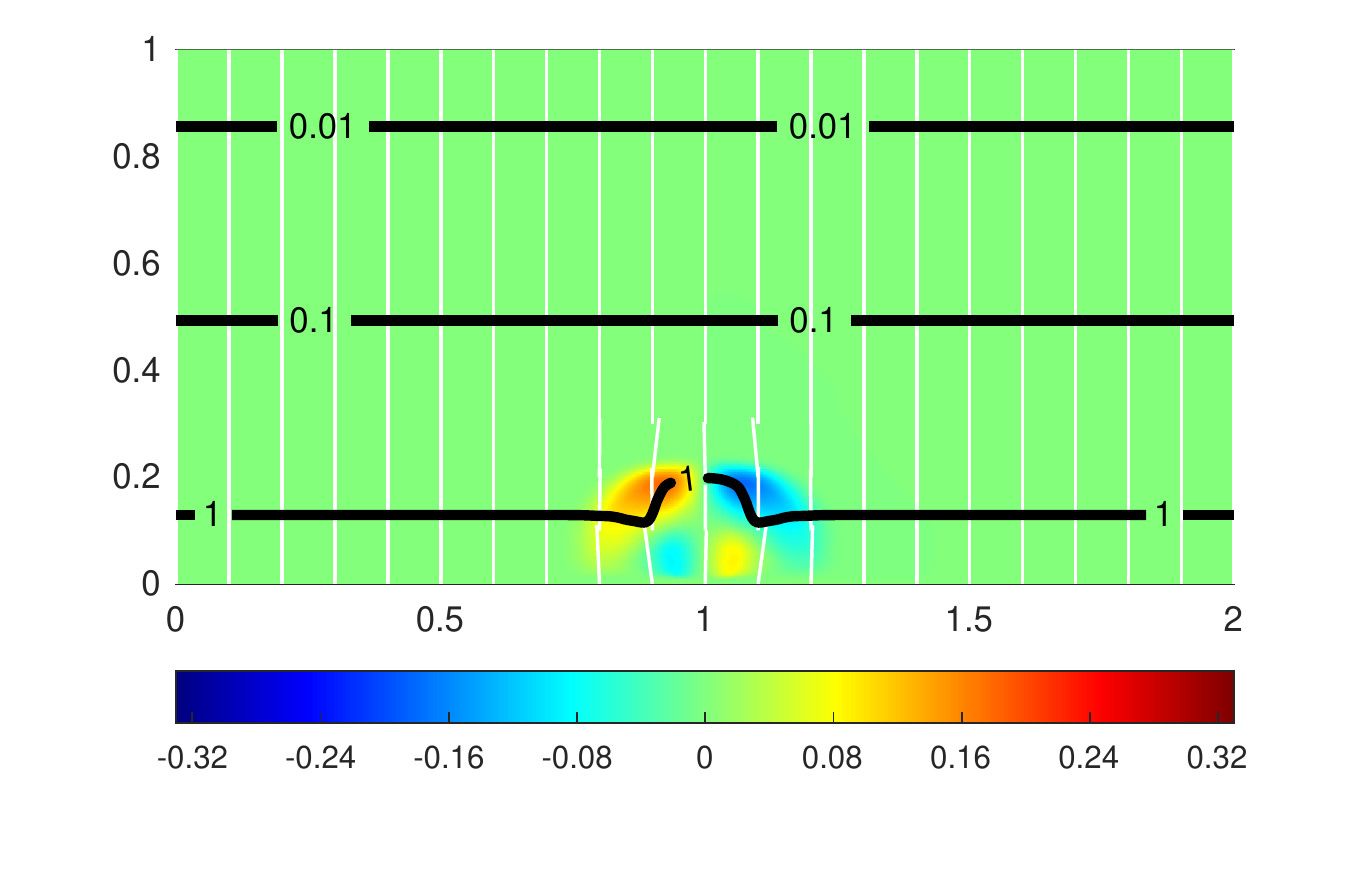}
		\subcaption*{$u_{\perp B}$ at $\text{t}=0.216$}
	\end{subfigure}%
	\begin{subfigure}[b]{0.3\linewidth}
		\centering
		\includegraphics[width=1\linewidth]{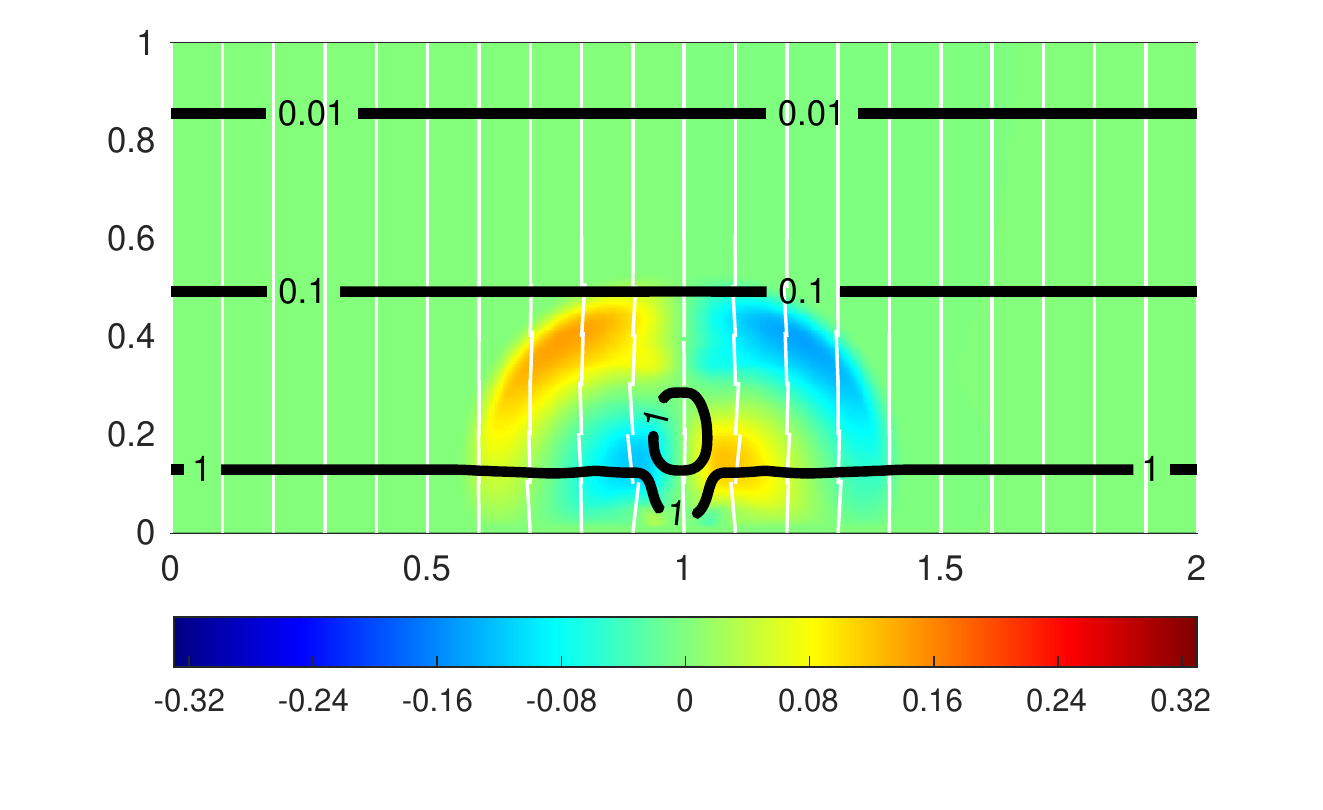}
		\subcaption*{$u_{\perp B}$ at $\text{t}=0.36$}
	\end{subfigure}%
	\begin{subfigure}[b]{0.3\linewidth}
		\centering
		\includegraphics[width=1\linewidth]{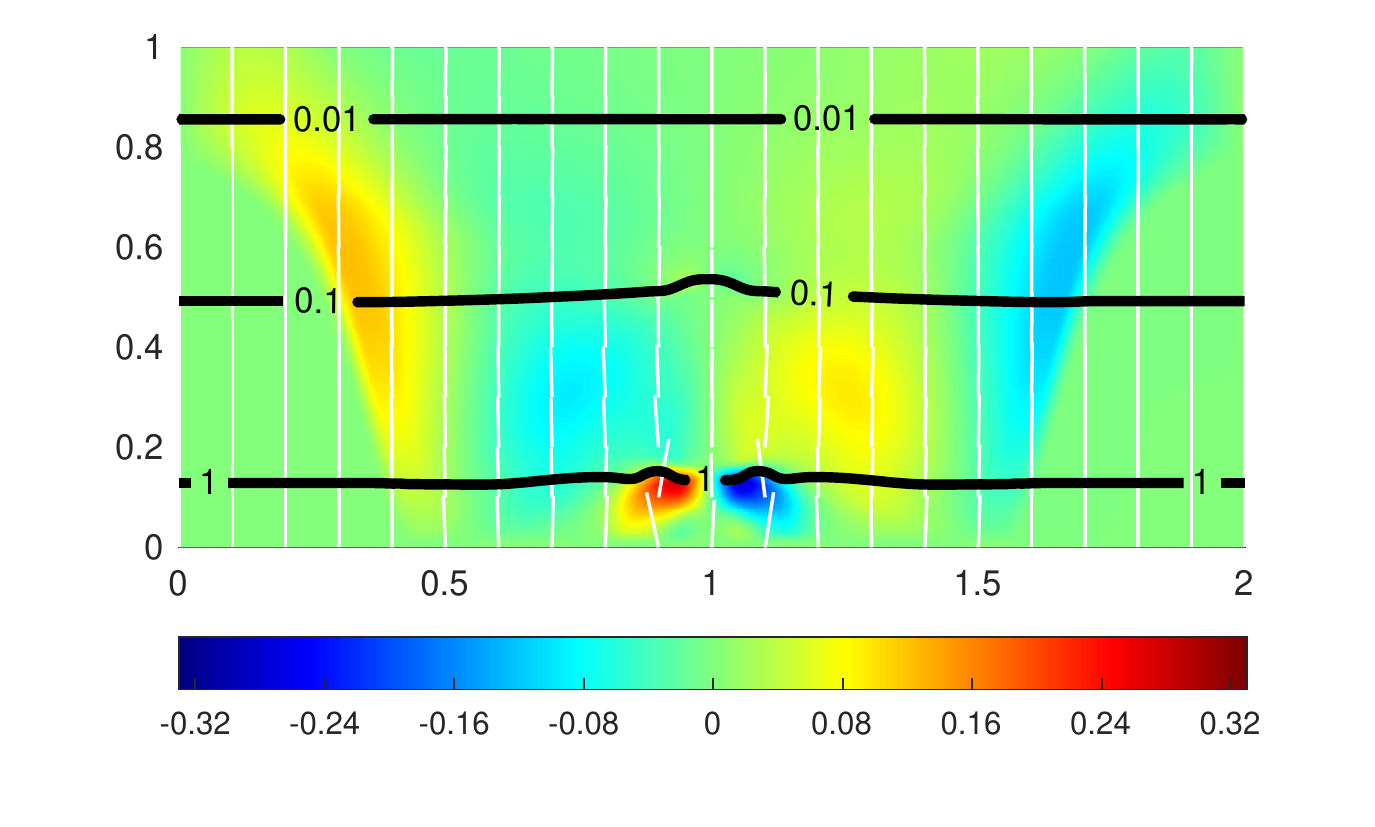}
		\subcaption*{$u_{\perp B}$ at $\text{t}=0.504$}
	\end{subfigure}
	\caption{MHD wave propagation: velocity perpendicular to the magnetic field $u_{\perp B}=<(u_1,u_2),(-B_2,B_1)>/|\textbf{B}|$ for $\mu=1$ on 400 $\times$ 200 grid points at different times.}
	\label{magnetoB}
\end{figure}

In this test case, we model propagating waves that not only undergo the effects of pressure and gravity, but also that of the magnetic field. The test case is extracted from \cite{Fuchs2010}. We consider the magnetohydrodynamic steady state defined as,\\
\begin{align}
	\rho(x,y)=\rho_{0}\exp(-\frac{y}{H}), p(x,y)=p_{0}\exp(-\frac{y}{H}),  \textbf{u}=0, \textbf{B}=(0, \mu,0), \nabla\cdot\textbf{B}=0.
\end{align}
Where $\mu$ is a parameter that takes different values for each part of the experiment. The waves model a perturbation of the steady state that starts from the bottom boundary of the normal velocity as follows,
\begin{equation}\label{piston}
	\textbf{u}_{i,\{0,1\}}^{n}=
	\begin{cases}
		\frac{\textbf{B}_{i,\{0,1\}}}{|\textbf{B}_{i,\{0,1\}}|}c\sin(6\pi t^n) \hspace{0.5 cm} \text{for } x \in [0.95, 1.05],\\
		0 \hspace{0.5 cm} \text{Otherwise,}
	\end{cases}
\end{equation}
with $c=0.3$. The computational domain is  $\left[0,2\right]\times\left[0,1\right]$. We use the wave propagation boundary conditions suggested in \cite{Fuchs2010}. These boundary conditions are periodic boundaries in the $x$-direction for $\textbf{U}$ and $p$ and Neumann type boundary conditions in the $y$-direction as the following, 
\begin{align*}
	\rho_{i,1}^n=\rho_{i,2}^ne^{\frac{\Delta y}{H}}, \rho_{i,0}^n=\rho_{i,1}^ne^{\frac{\Delta y}{H}}\\
	\rho_{i,ny-1}^n=\rho_{i,ny-2}^ne^{\frac{-\Delta y}{H}}, \rho_{i,ny}^n=\rho_{i,ny-1}^ne^{\frac{-\Delta y}{H}}
\end{align*}
for $1 \leq i\leq nx$. Similar boundary conditions for the momentum $\rho \textbf{u}$ and the pressure $p$. Energy boundary conditions are adopted from the pressure. For the magnetic field boundary conditions, we simply copy the data from the cell before. We present the profile of the velocity in the direction of the magnetic field,
\begin{align}
	u_B=<\textbf{u},\textbf{B}>/|\textbf{B}|,
\end{align}
at the final time $t=0.54$ for different values of $\mu$. As $\mu$ increases, the effect of the magnetic field on the propagating wave increases. The wave profile gets compressed as the magnetic field takes higher values. The plasma parameter is given by $\beta=\frac{2p}{\textbf{B}^2}$ \cite{Fuchs2010}. It measures the relative strength of the thermal pressure to the magnetic field, and is crucial in determining the dynamics of the plasma. The $\beta$-isolines are illustrated in black and the lines of the magnetic field are illustrated in white. The parameter $\beta$ indicates the effects of the pressure and the magnetic field on the propagating wave such that, for $\beta > 1$, the region is pressure dominated, while for $\beta < 1$, the region is magnetic field dominated. 
In figure \ref{magnetoB=0}, the profile of the velocity in the direction of the magnetic field, in the case of $\mu$ almost zero, is illustrated, which is exactly the velocity in the $y$-direction in this case. The wave propagates freely along the computational domain taking a radial profile in the absence of the magnetic field on 400 $\times$ 200 grid points. Figure \ref{magnetoB=1}, shows the profile of the propagating wave under the effect of a stronger magnetic field for $\mu=1$ on 400 $\times$ 200 grid points without applying CTM. In addition, figure \ref{magnetoB=1} presents the divergence of the magnetic field which is clearly not zero.  On the other hand, we present the same results with applying CTM on 1200 $\times$ 600 grid points in figure \ref{magnetoB=1CTM}. Applying the CTM results in a zero discrete divergence of the magnetic field up to machine precision. Another effect of applying the CTM is the diffusion we see in figure \ref{magnetoB=1CTM}, which was resolved by evolving the solution on a finer grid. 
Additionally, we present the  velocity in the direction perpendicular to the magnetic field in figure \ref{magnetoB} for $\mu=1$ at different times.\\
Our results, obtained with the second order scheme, are comparable with the results in \cite{Fuchs2010}, obtained with third order schemes, which ensures the robustness of our scheme and its capability of solving physically challenging problems, such as wave propagation under the effect of pressure and gravity. 

\section{Conclusion}
\label{sec::conclusion}
In conclusion, we develop  a two-dimensional second order unstaggered finite volume central scheme for the system of MHD equations. The proposed scheme is capable of preserving any type of known equilibrium states due to a special reformulation that computes the numerical solution in terms of a specific reference state. A comparison between the obtained numerical results and the corresponding literature ensures the robustness and the accuracy of the developed schemes. In this work, we chose the CTM as a procedure to clean the divergence of the magnetic field, which is applied dynamically whenever needed. Meaning that, in the test cases where the numerical divergence is zero at the final time and no numerical instabilities had been observed, we do not apply it. This leaves us with a second order well-balanced finite volume numerical scheme that captures solutions of the MHD equations and satisfies the divergence-free constraint.  All our computations are done on a Cartesian grid in 2D. A triangular mesh can be considered in future work.

\section{Funding}
\noindent The authors would like to acknowledge the National Council for Scientific Research of Lebanon (CNRS-L) for granting a doctoral fellowship to Farah Kanbar. Farah Kanbar also acknowledges funding by the Qualification Program of the Julius Maximilians University Würzburg.

\bibliographystyle{plain}







\end{document}